\documentclass[final]{siamltex}

\usepackage{amsmath}
\usepackage{amssymb}
\usepackage{graphicx}
\usepackage{epstopdf}
\usepackage{algorithm}
\usepackage{algorithmic,color}
\usepackage{multirow}
\usepackage{appendix}
\usepackage{arydshln}
\usepackage{cancel}
\usepackage{stackengine}
\usepackage{soul}

\newcommand{\R}{{\mathbb R}}

\newcommand{\N}{{\mathbb N}}

\newcommand{\tr}[1]{\textcolor{black}{#1}}

\newtheorem{defn}{Definition}

\newtheorem{thm}{Theorem}
\newtheorem{prop}{Proposition}
\title{Shape Partitioning via L$_\MakeLowercase{p}$ Compressed Modes}
\author{Martin Huska\thanks{
       Department of Mathematics, University of Bologna, Bologna, Italy.
        Email: martin.huska@unibo.it. }
        \and
        Damiana Lazzaro\thanks{
        Department of Mathematics, University of Bologna, Bologna, Italy.
        Email: damiana.lazzaro@unibo.it. }
        \and
				Serena Morigi\thanks{
        Department of Mathematics, University of Bologna, Bologna, Italy.
        Email: serena.morigi@unibo.it. }
        }

\date{Received: date / Accepted: date}
\begin{document}

\maketitle

\begin{abstract}
The eigenfunctions of the Laplace Beltrami operator (Manifold Harmonics) define a function basis
that can be used in spectral analysis on manifolds.
%of a Laplace-Beltrami operator can be realized as the solutions to variational optimization problems
In \cite{OLCO13} the authors recast the problem as an orthogonality constrained optimization problem and 
pioneer the use of an $L_1$ penalty term 
%to enforce compact support and 
to obtain sparse (localized) solutions.
In this context, the notion corresponding to sparsity is compact support which entails spatially localized solutions.
%is crucial, 
We propose to enforce such a compact support structure by a variational optimization formulation 
with an $L_p$ penalization term, with $0<p<1$.
%Extending previous results, we prove that for $0<p<1$, the eigenfunctions in the basis have compact support.
The challenging solution of the orthogonality constrained non-convex minimization problem is obtained by 
applying splitting strategies and an ADMM-based iterative algorithm.
The effectiveness of the novel compact support basis is demonstrated in the solution of the 2-manifold decomposition problem which plays an important role in shape geometry processing where the boundary of a 3D object is well represented by a polygonal mesh.
We propose an algorithm for mesh segmentation and patch-based partitioning (where a genus-0 surface patching is required).
Experiments on shape partitioning are conducted to validate the performance of the proposed compact support basis.
\end{abstract}

\bigskip

\section{Introduction}
\label{sec:intro}

Spectral analysis of the Laplace-Beltrami Operator (LBO) on a discrete manifold \tr{has} found 
many applications in surface processing, such as for example in shape matching, smoothing, shape 
recognition, and segmentation \cite{Reuter05,Zhang16,ZZ,WLT14,Song14,Mejia16}.

We are in particular interested in two significant surface processing applications 
of the LBO, namely mesh segmentation (or mesh partitioning) and patch-based partitioning,
which are members of a higher level class known as shape partitioning.

Mesh segmentation is fundamental for many computer graphics and animation
techniques such as modeling, rigging, shape-retrieval, and deformation. 
Given an object with arbitrary topology and a discrete manifold representing the object's boundary, 
this process consists in the decomposition of an object \tr{into} salient sub-parts and 
it relies mostly on surface geometric attributes of the object's boundary.

Patch-based partitioning has a variety of applications in product design and modeling, 
reverse engineering, texturing, and %rapid prototyping fields or
3D printing.
By means of this process a discrete manifold is decomposed into smaller patches or sub-manifolds
 that can be parameterized. Computations can then be performed on simple parameter domains.
The patch-based partitioning is often used together with a B-spline surface fitting technology.
In general patch-based patching produces discrete sub-manifolds of smaller size with respect to the mesh segmentation.

In recent years, the basis generated by the eigenfunctions of LBO, called Manifold Harmonic Basis (MHB)
has been proposed in \cite{VL08} in analogy to Fourier analysis, and used for example for object segmentation applications \cite{ZZ}.

However, in the shape partitioning context, rather than a multiresolution representation of the shape, which is the peculiarity of 
the MHB on manifolds, the focus is on identifying  the observable features of the manifold which represent 
for example protrusions, ridges, details in general localized in small regions.

Hence, in the partitioning context, a more suitable alternative to the MHB is represented by the Compressed Manifold Basis (CMB),
 introduced in \cite{NVTM}, which 
is characterized by compact support quasi-eigenfunctions of the LBO obtained by imposing sparsity constraints.

Motivated by the advantages in terms of control on the compact support obtained by using
the L$_1$ norm to force the sparsity of the solution discussed in \cite{OLCO13} and \cite{NVTM},
we devised to replace the $L_1$ norm by a more effective sparsity-inducing L$_p$ norm term, 
with $0 < p \leq 1$, which stronger enforces the locality of the resulting basis functions.
The set of functions $\Psi  = \left\{\psi_k\right\}_{k=1}^N$, that we will call $L_p$ Compressed Modes ($L_p$CMs),
 is computed  by solving the following variational model 
 \begin{equation}
%\left\{\psi_1,\psi_2,..,\psi_N\right\}=\arg 
\min _{\Psi}
\sum_{k=1}^{N} \left( \frac{1}{\mu} \|{\psi_k}\|_p^p -\frac{1}{2}<{\psi_k},\Delta {\psi_k}> \right) \quad s.t. \quad <{\psi_j},{\psi_k}>=\delta_{jk}
\label{problem}
\end{equation}
 where $\delta_{jk}$ is the Kronecker delta, and $\mu>0$ is a penalty parameter.
They form an orthonormal basis for the $L^2(\Omega)$ space, where $\Omega$ is the domain 
in consideration, 
%with respect to the standard inner product and 
and they represent a set of quasi-eigenfunctions of the Laplace-Beltrami operator.

%defines the $L_p$ Compressed Mode Basis.

%Concerning the regularization parameter $\mu \in \R$ in \eqref{problem}, it plays the crucial role 
%to control the compression of the support of the modes, such that
 %when the $\mu$ parameter increases the volume of the support
%of the $L_p$CMs shrinks, while decreasing $\mu$ produces small compression, i.e. the enlargement of the support.

The second term in the objective function of \eqref{problem} is the fidelity term which represents the accuracy of the shape approximation provided by the set of 
functions $\Psi$, while the first term, so-called penalty term, 
forces the sparsity  in the functions $\Psi$ thus imposing spatially sparse solutions.
We remark that at the aim to construct a basis which is sparse but also localized in space it is necessary 
to further demonstrate that the functions $\Psi$ determined by solving \eqref{problem} have compact support.
This aspect will be proved in this work.

The penalty parameter $\mu$ controls the compromise between the two aspects.
It is well known that the sparsity is better induced by the L$_p$ norm for $0<p<1$, rather 
than the L$_1$ norm.
For $p=1$ model \eqref{problem} reduces to the proposal in \cite{OLCO13}, 
where the sparsity is forced only by acting on the $\mu$ value
to increase the contribution of the penalty term, thus decreasing the shape approximation
guaranteed by the fidelity term.
 
The parameter $p$ plays a crucial role since it allows to force the sparsity while maintaining 
the approximation accuracy without excessively stressing the penalty via the $\mu$ value. 
The accuracy is fundamental to localize the support of the functions in specific local features of the shape
such as protrusions and ridges.   		
\begin{figure}
	\flushleft
	\begin{minipage}{\textwidth}
		\centering
		\includegraphics[width=2.0cm]{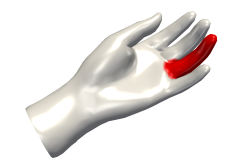} 
		\includegraphics[width=2.0cm]{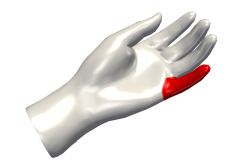} 
		\includegraphics[width=2.0cm]{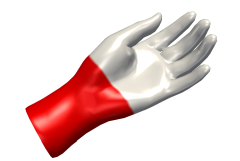} 
		\includegraphics[width=2.0cm]{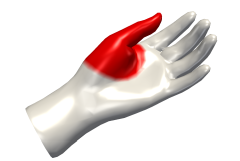} 
		\includegraphics[width=2.0cm]{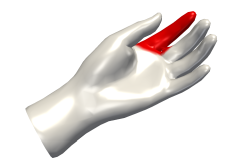} 
		\includegraphics[width=2.0cm]{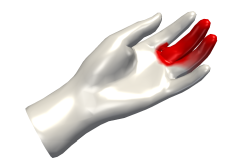}   
		
		(a) Supports of L$_p$CM Basis, $p$=1, $\mu=100$;
	\end{minipage}
	\begin{minipage}{\textwidth}
		\centering
		\includegraphics[width=2cm]{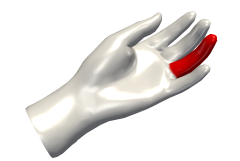} 
		\includegraphics[width=2cm]{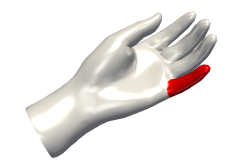} 
		\includegraphics[width=2cm]{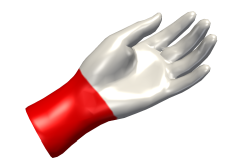} 
		\includegraphics[width=2cm]{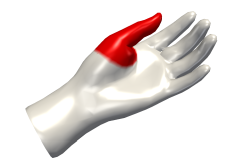} 
		\includegraphics[width=2cm]{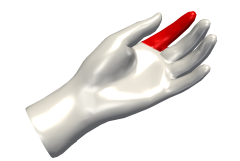} 
		\includegraphics[width=2cm]{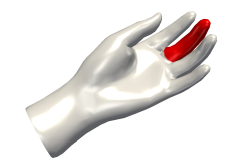} 
		
		(b) Supports of L$_p$CM Basis, $p=0.8$, $\mu=1000$;
	\end{minipage}			
	
	\caption{Partitioning of the 2-manifold \texttt{hand} using L$_1$CM basis (a) and L$_p$CM basis (b).}
	\label{fig:hand}
\end{figure}
	
Some evidence of the benefit obtained by the sparsity-inducing proposal, is shown in Fig.\ref{fig:hand}
	where we try to answer the following question: Can we identify the most salient parts of the manifold \texttt{hand} using only six compressed modes?
Fig.\ref{fig:hand} compares the supports of the compressed modes determined as solution of  the variational problem \eqref{problem} with $p=1$ (Fig.\ref{fig:hand}(a)) and $p=0.8$ (Fig.\ref{fig:hand} (b)) 
where the $\mu$ parameter has been \tr{chosen, for each $p$ value, to provide the most natural salient part identification.}
An automatic strategy for the choice of the optimal $\mu$ parameter will be also discussed in this work.
The supports of the  six quasi-eigenfunctions are colored in red and we can observe 
	that using $p<1$ strengthens the sparsity, while if $p=1$ 
no $\mu$ value has allowed to correctly identify all the fingers.
%allowing for a better localization on the shape features. 

\smallskip

An efficient solution of the orthogonality constrained problem (\ref{problem}) represents a challenging task
due to the nonlinear, non-convex orthogonality constraints combined with the non-smooth and non-convex objective function,  
which may lead to many different local minimizers as solutions. 
Non-trivial iterative approaches are commonly used to solve this kind of optimization problems.
In this paper we propose a variant of the basic Alternating Direction Method of Multipliers  (ADMM) approach \cite{BOYD_ADMM}
where the non-convex orthogonality constraints defined for the $L_p$ Compressed Modes $\Psi$ in (\ref{problem})
are preserved by means of a SVD matrix factorization, following \cite{LO14},
and a suitable proximal operator is devised to deal with the non-convex penalty.

Although we currently have no proof of convergence for this method, our numerical experiments 
verify that the method generates a sequence of iterative solutions preserving the orthogonality constraints and decreasing the cost functional 
in (\ref{problem}).  This anyway will motivate future theoretical analysis.

\smallskip

Finally, the proposed ADMM algorithm becomes the kernel of a more sophisticated strategy 
finalized to the decomposition of a 2-manifold approximated by only a few L$_p$CMs.
The \tr{identification and selection} of features is a key problem in the context of shape partitioning 
and the L$_p$CMs hold the potential for naturally handle this problem thanks to the compact support 
property that characterizes them.

We devised two methods for generating a suitable L$_p$CM basis, the first 
imposes the $\mu$ parameter in \eqref{problem} while leaving the number of functions free to increase
until the shape boundary is completely covered; the second method is based instead on an automatic tuning of the parameter $\mu$ 
for a fixed number of L$_p$ compressed modes. Once the set of functions $\Psi$ is determined, 
a region growing process is applied to construct a partitioning of the mesh.
\tr{A unified algorithm is proposed for both mesh partitioning and a more demanding
patch-based partitioning (which guarantee genus-0 patches). }

%The proposed algorithms run very fast and return solutions no worse than those from their state-
%of-the-art algorithms.
 	
	%
%\begin{enumerate}
  %\item 
		%
	%In many application situations, a global parameterization that maps a surface onto one
%global parameter domain is desirable. The global nature avoids the needs of introducing cuts
%on the surface and partitioning the surface into several patches. Computations can then be
%performed on one simple parameter domain.
	%basis use of (\ref{problem}) for segmentation and patch-based partitioning
	  %
%\end{enumerate}

Summarizing, the main contributions are as follows:
\begin{itemize}
\item[a)] a new variational model for the construction of a compact support basis for the Laplacian operator;
\item[b)] analysis of the compact support property of the obtained basis functions;
\item[c)] proposal of an efficient algorithm for the solution of the optimization problem 
based on ADMM; 
\item[d)] devise of a partitioning algorithm for both shape segmentation and patch-based partitioning based on L$_p$ compressed modes.
\end{itemize}
 
\bigskip

The paper is organized as follows. In Section \ref{sec:cm} we briefly review 
the compressed modes and their extension to the 2-manifold context. 
In Section \ref{sec:prop} we introduce the sparsity-inducing variational model  
to determine the L$_p$CMs and we provide two important theoretical results.
The discretization of the optimization problem on triangulated surfaces is 
given in Section \ref{sec:discr}. 
An efficient ADMM-based iterative algorithm for the solution of the discrete version of (\ref{problem}) is presented in Section 
\ref{sec:ADMM}. Basic notions on the partitioning problem are provided in Section \ref{sec:basic}
and the algorithmic proposal for both segmentation and patch-based partitioning is described 
in Section \ref{sec:alg}.
Numerical experiments are presented in Section \ref{sec:ne} and 
conclusions are drawn in Section \ref{sec:conc}.

\subsection{Related work on shape partitioning}
Patch-based partitioning and mesh segmentation have been widely studied over recent years,
	creating a whole categorization of methods based on different methodologies, see  \cite{MSsurvey,CADA,Chen,AKM,Sha08}.
%The idea in the mesh segmentation context localize the meeting of two salient parts of the object in its concave areas.
%Such segmentation methods based on concavity-awareness is proposed in \cite{AZC*12} and \cite{mLA07}. 
%In Asafi et al. \cite{AGCO13} weakly convex components are obtained by a point-visibility test.

In order to partition a mesh, the spectral analysis methods use the eigenvalues of properly defined matrices,
called affinity matrices,  based on the connectivity of the mesh. 
The authors in \cite{LZ04} define an affinity matrix using both geodesic and angular distances.
The spectral analysis is performed on the Laplacian matrix weighted by 
 dihedral angle differences in \cite{ZZ}, and by mean curvature in \cite{HM},
 and then successfully applied to	mesh partitioning.
\tr{In \cite{Chau} an affinity matrix is proposed based on the optimal normalized Cheeger cut
which encodes both the structural and the geometrical information in order to segment concave regions.}

The results of the methods in the class of spectral-based mesh segmentation strongly depend on
		the affinity matrix considered. In the proposed approach good quality results are obtained simply
		from the Laplace-Beltrami spectral decomposition, by imposing suitable constraints of orthogonality and sparsity.

%Variational partitioning has been introduced in \cite{CSAD04} where the authors presented an optimization cost function based on clustering face normal of the mesh. 
%Several variational-based partitioning methods have been proposed since then to solve both patch-based partitioning
%and mesh segmentation. 

In \cite{MSsurvey}, \textit{patch-based partitioning} is named surface-type segmentation and 
what is here defined as \textit{mesh segmentation} is instead referred as part-type segmentation.

An important result on part-type segmentation has been presented in \cite{ZZ}, where a convexified
version of the variational Mumford-Shah model is presented and extended to 3D meshes.
In reverse engineering and Computer Aided Design (CAD) applications, patch-based partitioning is seen as an
 automatic procedure to create CAD models from measured data, \cite{varady}.
Patch-based partitioning in \cite{YWLY} is finalized to fitting quadric surfaces to the mesh,
while in \cite{FIT} a smooth stitch of spline patches is built. 
However, the success of any patch-based method strongly depends on the goodness of the underlying shape partitioning.		

\tr{We propose a unified framework to perform both mesh segmentation and patch-based partitioning.
The proposed method carries out two approaches, one completely unsupervised, namely the number of segments is determined automatically, 
and the other supervised, by performing a segmentation into a given number of parts. 
We refer the reader to \cite{The} for unsupervised state-of-the-art methods,  and to \cite{Ben,Kal} for supervised competitors.}
\section{Background on Compressed Modes}
\label{sec:cm}

In the preliminary work \cite{OLCO13} the authors show how to produce 
a basis of localized functions $\{\psi_k\}_{k=1}^N$ in $\R^d$, called Compressed Modes (CMs),  
by solving the following variational problem 
\begin{equation}
 \min_{\left\{\psi_1,\psi_2,..,\psi_N\right\}}
 \sum_{k=1}^{N} \left( \frac{1}{\mu} \|{\psi_k}\|_1 + \left\langle {\psi_k},H {\psi_k} \right\rangle \right) 
\quad s.t. \quad <{\psi_j},{\psi_k}>=\delta_{jk},
\label{CM}
\end{equation}
where $H= -\frac{1}{2} \Delta + V(x)$ is the Hamiltonian operator corresponding to potential $V(x)$,
 the $L_1$ norm is defined as $\|f\|_1 = \int_{\Omega} |f| dx$ and $\left\langle f,g \right\rangle = \int_{\Omega} f^*g dx \; \Omega  \subset \R^d$. Here the L$_1$ norm is a penalty term used to achieve spatial sparsity.
The orthonormality constraints in (\ref{CM}), which enforce  the orthonormality of the basis functions,
lead to a non-convex variational problem, with many local minimizers. 
%a space of feasible functions which is not a convex set,   

A theoretical analysis of the CMs, provided in \cite{Osher}, allows for 
finding the minimizer of the variational formulation of the Schrodinger equation, 
showing the spatial localization property of CMs, and
establishing an upper bound on the volume of their support.
Consistency results for the CMs %constructed from the variational formulation
%of the Schrodinger equation 
were proved in \cite{F14}.

In \cite{NVTM} the variational problem (\ref{CM}) in $\R^d$ domains has been extended 
to deal with Laplace-Beltrami eigenfunctions on 2-manifolds discretized by three-dimensional meshes.
These new basis functions, named Compressed Manifold Modes (CMM), 
form the Compressed Manifold Basis (CMB) and define an alternative to the classical MHB, \tr{proposed} in \cite{VL08}.
It is well known that  the eigenfunctions of the Laplace Beltrami operator, called Manifold Harmonics (MH), 
define a function basis. In particular,
for a smooth manifold $\mathcal{M}$ embedded in $\R^3$ the
Laplace-Beltrami operator induces a set of eigenfunctions $\{ \phi_k\}$ and associated eigenvalues 
$\{\lambda_k\}$ 
determined by 
\begin{equation}
-\Delta \phi_k = \lambda_k \phi_k \qquad k \in \N , \; \lambda_k \in \R.
\label{MHB}
\end{equation}
The self-adjointness of $\Delta$ implies that the eigenvalues are real and that the eigenfunctions are orthogonal 
with respect to the L$_2$-inner product: $<f,g>=\int_{\mathcal{M}} f \; g$.   
%These functions  on the manifold form a basis, which has been named MHB in \cite{VL08}, where an efficient numerical
%algorithm is proposed to compute the MHB on a triangular mesh representing the discretization of the manifold $\mathcal{M}$.

One major drawback of this basis is that, similarly to the Fourier spectrum, the MHs are dense and have global spatial support. 
This means that the functions do not give intuitive insight on the features of the manifold, thus reducing their
use in practical shape processing applications \cite{LZ10}. 
However, it is well known that using a reduced number of eigenfunctions corresponding to the smallest eigenvalues $\lambda$,
the MHs allow to approximate the shape of the manifold in an improved manner as the number of eigenfunctions increases.

%
%\bigskip
%%
%\textit{The matrix $\mathcal{L} \in \R^{n \times n}$ defined in (\ref{eq:lap}) associated to a connected mesh $\Omega$ of $n$ vertices,
%satisfies the following properties:
%%
%\begin{description}
%\item{ 1)} $\mathcal{L}$ is symmetric and positive semi-definite;
%\item{ 2)} $\mathcal{L}=U\,\Lambda\,U^T, \, \Lambda = diag(\lambda_i), \, 0\,= \lambda_0 < \lambda_1 < \ldots <\lambda_n$;
%\item{ 3)} $\lambda_i, \forall i$ are real eigenvalues, $U^T\,U\,=\,I_n$ with $I_n$ the identity matrix of order $n$, $U=\{v_0,v_1,\ldots,v_n\}$ 
%form an orthogonal basis of $\R^n$;
%\item{ 4)} If $f\,=\, \sum_{i=1}^n \left\langle f, v_i \right\rangle v_i$, the $k$-term approximation of $f$ is given by
%$$f_k\,=\, \sum_{i=1}^k \left\langle f, v_i \right\rangle v_i.$$   
%\label{eq:eigen_vect}
%\end{description}}
%The spectral decomposition of $\mathcal{L}$, defined in the following, % Proposition \ref{prop:eigen},
%provides a set of $(n-1)$ non trivial, smooth, shape intrinsic isometric-invariant maps. 
%We refer the reader to \cite{VL} for details.

\section{The sparsity-inducing variational model for L$_p$CMs}
\label{sec:prop}

Let $\Omega=B(0,R) \subset \R^d$ \tr{denote} the $d$-dimensional ball of radius $R$ centered at the origin.
Relation (\ref{problem}) can be rewritten as:
\begin{equation}
\label{problem1}
%\left\{\psi_1,\psi_2,..,\psi_N\right\}=\arg 
\min_{\Psi}
\sum_{k=1}^{N} \left( \frac{1}{\mu} \int_\Omega{\left|{\psi_k}\right|^p dx} -\frac{1}{2} 
 \int_\Omega{{\psi_k} \Delta {\psi_k} dx}\right) \quad s.t. \quad \int_\Omega{{\psi_j} \, {\psi_k} \,dx} \, = \, \delta_{jk},
\end{equation}
where we denoted the $L_p$ norm of a function by $\| f\|_p = \left( \int_{\Omega} |f|^p dx \right)^{1/p}$.

%\begin{thm} There exist a minimizer to problem (\ref{problem2}).
%\end{thm}
%\subsection{Bounds on the Volume of Support of Compressed Modes}

Before proceeding with the solution of the variational model \eqref{problem1} we demonstrate
 in the following two important properties of the L$_p$ compressed functions such as local support and completeness,
which also hold for the CMs determined by solving  \eqref{CM}.

\subsection{On the support of the L$_p$CMs}

In this subsection we establish an asymptotic upper bound on the volume of the support of the L$_p$CMs in terms of the penalty 
parameter $\mu$ and the sparsity parameter $p$.
At this aim, we first reformulate (\ref{problem1}) by using integration by parts and imposing zero boundary conditions on $\Omega$.
It follows that the first $N$ compressed modes $\left\{\psi_i \right\}_{i=1}^N$ solve the following constrained optimization problem:
\begin{equation}
\label{problem2}
 \min_\Psi
\sum_{i=1}^{N} \left( \frac{1}{\mu} \int_\Omega{\left|{\psi_i}\right|^p dx} +\frac{1}{2} 
 \int_\Omega{\left| \nabla {\psi_i} \right |^2 dx}\right) \quad s.t. \quad \int_\Omega{{\psi_j}{\psi_k}dx}=\delta_{jk}.
\end{equation}

We first introduce the following result on the volume support of the first compressed mode. 

\smallskip

\begin{prop}
\label{extProp3_4}
$\forall x \in \Omega$, any $0 <p$ \tr{ $<$}  $1$,  and $\mu$ sufficiently small, we have 
\begin{equation}
\label{volume}
\int_\Omega{\left( \frac{1}{\mu}\left|{\psi_1}\right|^p+\frac{1}{2}\left |\nabla {\psi_1} \right |^2 \right)  dx}\leq m(\Omega)^{\frac{1}{p}-1}\mu^{-\frac{4}{4+d}}
\end{equation}
where $m(\Omega)$ is the finite measure of the domain $\Omega \subset \R^d$.
\end{prop}

\smallskip

\begin{proof}
From the relation between $L_p$ norm and $L_q$ norm, with $0<p<q\leq \infty$
\begin{equation}
\label{lplq}
\|f\|_p \leq m(\Omega)^{\frac{1}{p}-\frac{1}{q}}\|f\|_q,
\end{equation}
for $q=1$ and $0<p<1$, it follows that 
\begin{equation}
\int_\Omega{\left( \frac{1}{\mu}\left|{\psi_1}\right|^p+\frac{1}{2}\left|\nabla {\psi_1} \right |^2 dx \right)} 
\leq m(\Omega)^{\frac{1}{p}-1} \cdot\int_\Omega{\left( \frac{1}{\mu}\left|{\psi_1}\right|+\frac{1}{2}\left |\nabla {\psi_1} \right |^2 \right) dx.}
\end{equation}
By using Proposition 3.4 of \cite{Osher}, namely
$$\int_\Omega{\left( \frac{1}{\mu}\left|{\psi_1}\right|+\frac{1}{2}\left|\nabla {\psi_1} \right |^2 \right)dx}=C_1 \mu^{- \frac{4}{4+d}},$$
where $C_1$ is some fixed constant depending on $d$, we easily obtain the bound in \eqref{volume}.
%can conclude:
%\begin{equation}
%\int_\Omega{\left( \frac{1}{\mu}\left|{\psi_1}\right|^p+\frac{1}{2}\left|\nabla {\psi_1} \right |^2 dx \right)} \leq m(\Omega)^{\frac{1}{p}-1} \cdot C_1 \cdot \mu^{-\frac{4}{4+d}}
%\end{equation}
\end{proof}

\bigskip

\begin{thm}
There exist $\mu_0$, such that for $\mu < \mu_0$ the corresponding L$_p$ compressed modes $\left\{\psi\tr{_i} \right\}_{i=1}^N$ satisfy
\begin{equation}
\left | supp (\psi_i) \right|  \leq C \mu^{-\frac{8}{4+d}+1} m(\Omega) ^{\frac{1}{p(1-p)}-2}
\end{equation}
\label{th:cs}
where $C$ depends on $N$ and $p$.
\end{thm}

\bigskip

The proof is postponed to the appendix.

\bigskip

The result in Theorem \ref{th:cs} is fundamental for the construction of a compact support L$_p$CM basis  
and it will represent the key aspect for the shape partitioning method based on the L$_p$CMs
that will be described in Section \ref{sec:alg}.

\begin{figure*}[ht]
	\centering
	\includegraphics[width=3.5cm]{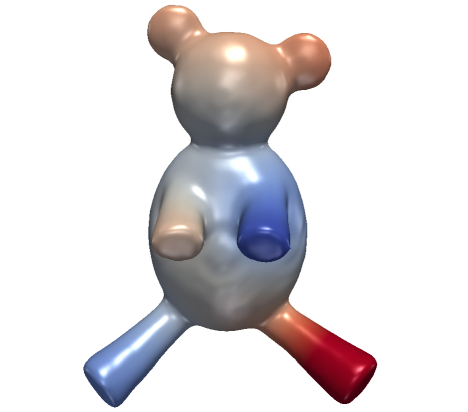} \hspace{-0.7cm}
	\includegraphics[width=3.5cm]{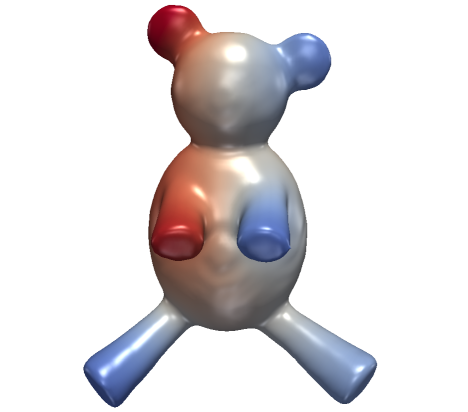} \hspace{-0.7cm}
	\includegraphics[width=3.5cm]{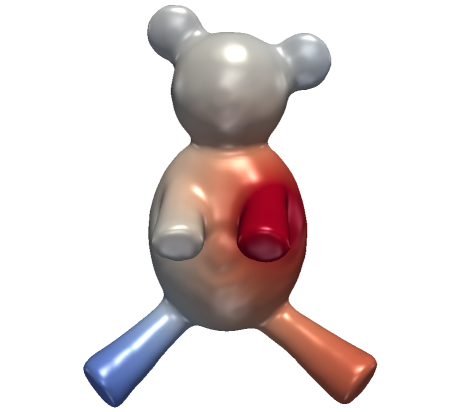} \\
	\includegraphics[width=3.5cm]{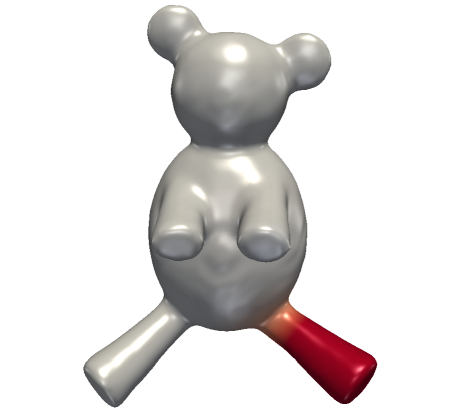} \hspace{-0.7cm}
	\includegraphics[width=3.5cm]{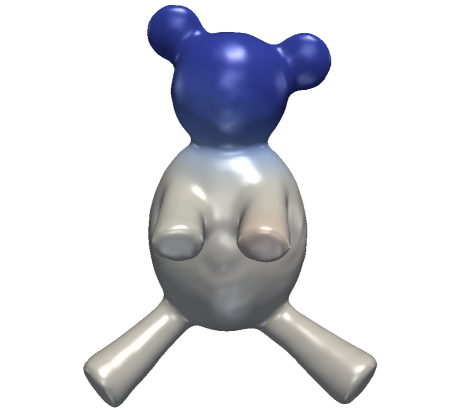} \hspace{-0.7cm}
	\includegraphics[width=3.5cm]{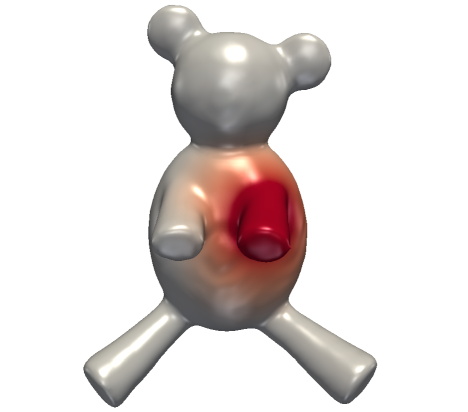} \\
	\includegraphics[width=3.5cm]{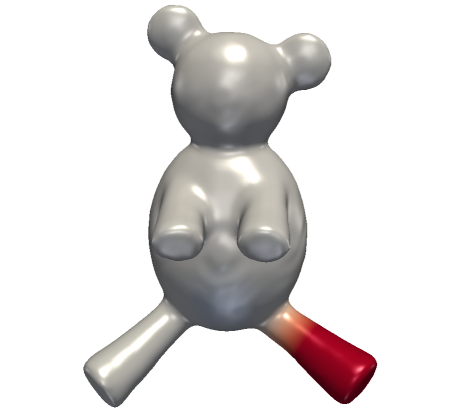} \hspace{-0.7cm}
	\includegraphics[width=3.5cm]{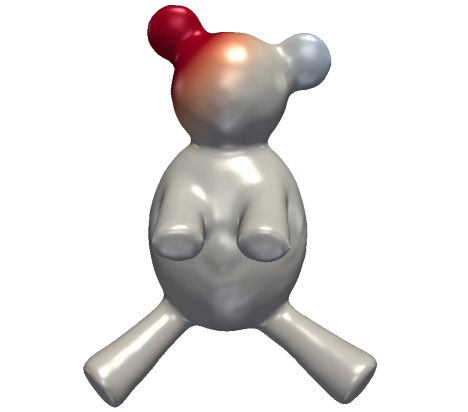} \hspace{-0.7cm}
	\includegraphics[width=3.5cm]{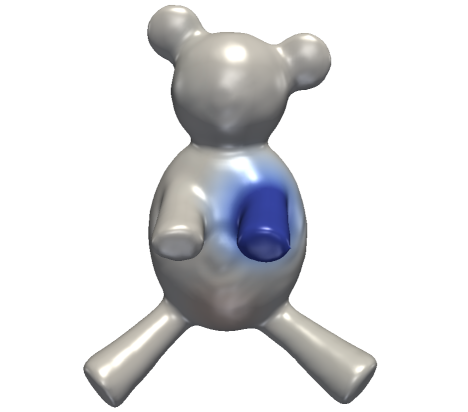} 
	
	\caption{ Three L$_p$CMs generated for the \texttt{teddy\_2} mesh for three different $\mu$ values,
		row-wise $\mu=\{250,160,50\}$ }
	\label{fig:Teddy}
\end{figure*}

\tr{An example demonstrating the essence of Theorem \ref{th:cs} is shown in Figure \ref{fig:Teddy}. 
In each row three L$_p$CM functions are illustrated obtained for a particular $\mu $ value, and fixed  $N$ and $p$ values.
Since the upper bound given in \eqref{th:cs} depends on $\mu, p, N$, for increasing values of $\mu$, as we expected, 
we notice an enlargement of the compact support of each function.
}

%\phantom{aaa}\newline \\
\subsection{Completeness of the $L_p$CMs} %\\

We now investigate a completeness result on the $L_p$CMs and its effect on shape approximation. 
In particular, we prove that, for a fixed $\mu$ value in (\ref{problem1}), under some unitary transformations, 
the $L_p$CM functions $\left\{ \psi_i\right\}_{i=1}^N$ approximate the eigenfunctions of the Laplacian operator 
in an improved manner as $N$ increases.
%\\

Let $\Phi=\left\{ \phi_i\right\}_{i=1}^M$ be the set of orthonormal eigenfunctions of $-\frac{1}{2}\Delta$ corresponding to the eigenvalues $\left\{\lambda_i\right\}_{i=1}^M$, defined by (\ref{MHB}) where the eigenvalues are arranged in non-decreasing order.
%\\

The following result of completeness for the $L_p$CMs holds.
%\\

\smallskip

\begin{thm}
Given fixed parameters $\mu$ and $p$ in (\ref{problem1}), for a fixed integer $M<N$, the first $N$ functions $L_p$CMs $\left\{\psi_i \right\}_{i=1}^N$ up to an 
unitary transformation,  satisfy
\begin{equation}
\lim_{N \rightarrow \infty} \|\phi_i-\psi_i\|_2^2 =0, \quad i=1,\ldots,M.
\end{equation}
\label{th:completeness}
\end{thm}

\smallskip

\begin{proof}
The proof follows from \cite{YO13}, where the authors demonstrate the result in the case of $L_1$-norm, 
but it still holds if the $L_1$ norm term is replaced by any functional bounded by
$L_2$ norm. In fact, for the relation \eqref{lplq} between $L_p$ and $L_q$ norms, $0<p<q\leq \infty$, 
if we set $q=2$, it follows that $L_p$-norm, $0<p<2$, is bounded by $L_2$-norm.
\end{proof}

The completeness result confirms \tr{that} using the $L_p$CM orthogonal basis, analogously to the $\Phi$ basis,
we can reconstruct any function defined on the shape, 
up to an arbitrary degree of precision.   
However, for a small number $N$ of functions, the approximated reconstructions show significant differences. 

\begin{figure}
	\centering
	\includegraphics[width=4.1cm]{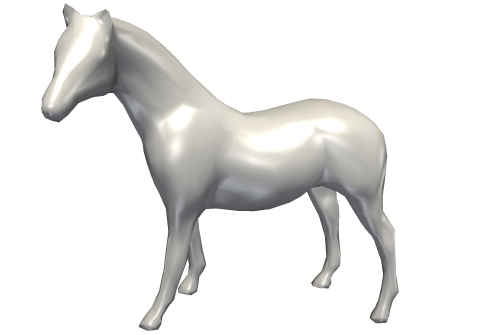} \\
	\includegraphics[width=4.1cm]{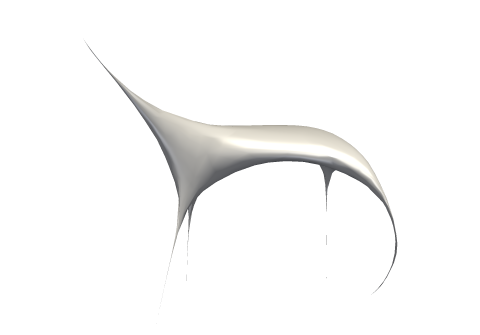}	 %\hspace{-0.5cm}
	\includegraphics[width=4.1cm]{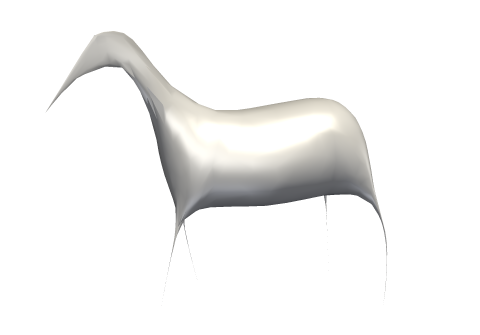}	 %\hspace{-0.5cm}
	\includegraphics[width=4.1cm]{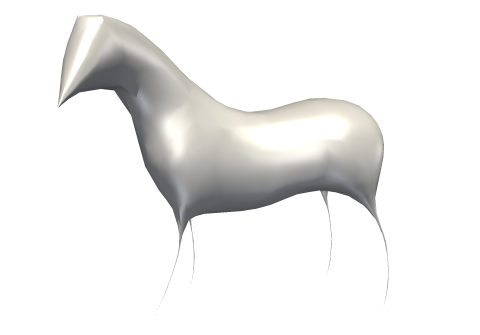}	\\ %\hspace{-0.5cm}
	\includegraphics[width=4.1cm]{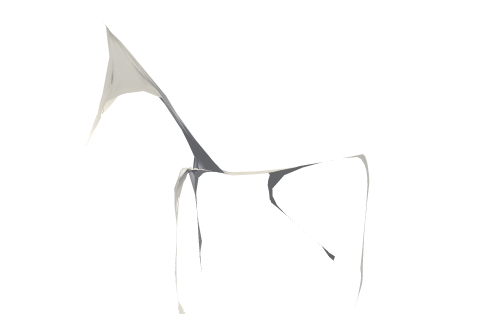}  		
	\includegraphics[width=4.1cm]{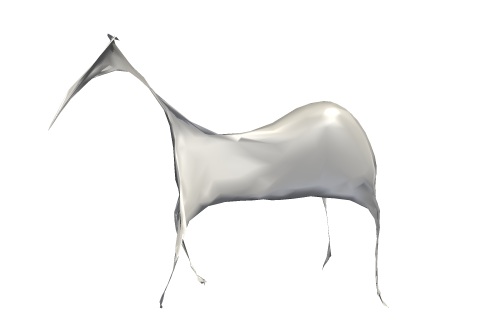} 		
	\includegraphics[width=4.1cm]{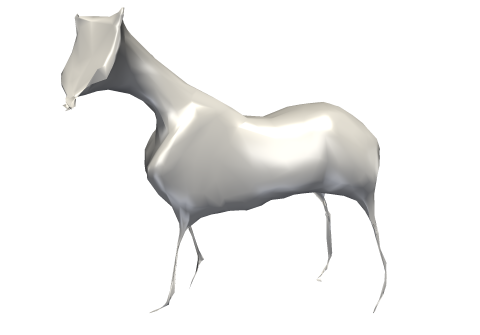} 		

	\caption{Reconstruction of the \texttt{horse} shape (top) using MH (middle row) and L$_p$CM (bottom row) bases 
	for increasing space dimension $N$; from left to right $N=8$, $N=15$, and $N=30$. }
	\label{fig:recon}
\end{figure}

By the way of illustration, let us consider the geometric reconstruction %of %the shape represented by the LBO
 of the 2-manifold \texttt{horse}.% illustrated in Fig.\ref{fig:recon} (top).
The shape approximation process, described for MHs in \cite{VL08}, also holds for L$_p$CMs.
The reconstruction obtained by using all the $N$ eigenfunctions of the LBO, where $N=868$,
 is shown in Fig.\ref{fig:recon} (top).
In the second and third row of Fig.\ref{fig:recon} we show, respectively,  the shape reconstructions obtained using the basis 
MHB computed by solving  \eqref{MHB},  and the proposed L$_p$CM basis obtained by \eqref{problem1},
formulated in the 2-manifold context which will be \tr{discussed} in the following sections. 
For a fixed value of the dimension $N$ in the range $N=8,15,30$, 
the reconstruction obtained by MHB is smoother but less representative of the \tr{underlying} shape, 
while the  L$_p$CM approximation \tr{looks like a more stylized shape, roughly a skeleton of the shape}. 
Moreover, while the MHB approximations for increasing dimensions tend to refine the basic shape, the L$_p$CM basis  
enriches the skeleton shape with smaller features while maintaining the structure of the shape.
This can be observed in the horse reconstructions in Fig.\ref{fig:recon} (bottom) where the 
horse's legs and ears are well represented only using the L$_p$CM bases.

%\\

%
%\begin{figure}
		%\centering
			%\begin{tabular}{cccc}
	%& \includegraphics[width=3cm]{figs/CM_1_1.png}	& %\hspace{-0.5cm}
	%\includegraphics[width=3cm]{figs/CM_1_2.png}	& \\ %\hspace{-0.5cm}
	%\includegraphics[width=3cm]{figs/CM_1_3.png}&	 %\hspace{-0.5cm}
	%\includegraphics[width=3cm]{figs/CM_1_4.png} & 		
	%\includegraphics[width=3cm]{figs/CM_1_5.png}& 		
	%\includegraphics[width=3cm]{figs/CM_1_6.png} \\		
		%\end{tabular}
		%\caption{Support of L$_1$ CM Bases $\mu$=0.0001}
		%%\label{fig:recon}
	%\end{figure}
	%

%\\
%\begin{thm}
%As $\mu \rightarrow \infty $, the eigenvalues $\sigma_i$'s of the matrix
%$<\Psi_N, \Delta \Psi_N>$, whose $(j,k)-th$ entry defined by
%$\int_\Omega{\psi_j \Delta \psi_k}$ converge to the eigenvalues $\lambda_i$ of $-\frac{1}{2}\Delta.$
%\label{th:completeness1}
%\end{thm}

\section{Discretization of the variational model}
\label{sec:discr}

We are interested in the application of model (\ref{problem}) to compute the $L_p$ Compressed Modes 
induced by the LBO on a manifold $\mathcal{M}$. We approximate $\mathcal{M}$ by a triangulated surface 
mesh $M := (V,T)$, where $V=\{X_1,\cdots,X_n \}$ is the set of $n$ vertices, $T$ is the connectivity graph,
and we denote by  $\mathcal{E} \subseteq {V} \times {V}$ the set of edges.
Each vertex $X_i \in V$ has immediate neighbors $X_j, j \in N(X_i),$ to which it is connected by a single edge $e_{ij}$.
We denote by $ N_\triangle(X_i)$ the set of triangles with vertex $X_i$, and by 
$| N_\triangle(X_i) | := \sum_{j\in N_\triangle(X_i)}A(\tau\tr{_j})$,  where $A(\tau\tr{_j})$ is the area of the triangle 
$\tau\tr{_j}$. 

We first introduce a popular discretization of the Laplace-Beltrami operator for a triangle mesh $M$, 
which, according to \cite{PP}, may be realized by $D^{-1}L$, where $L \in \R^{n\times n}$ 
is a symmetric, positive semi-definite, sparse matrix (weight matrix) defined as 
\begin{equation}
\begin{array}{l}
		L(i,j) := \left\{ 
		\begin{array}{ll}  
		\omega_{ij} = \frac{1}{2}(\cot\gamma_j + \cot\delta_j) & j \in N(X_i) \\
		-\sum_{k \in N(X_i)} \omega_{ik} & i=j \\
		0 & otherwise
		\end{array}\right.
\label{eq:N-mesh_1}
\end{array}	
\end{equation}	
%\begin{equation}
%\begin{array}{l}
		%L(X_i) = \dfrac{1}{2}\sum\limits_{j\in N(X_i)}\omega_{ij}\left(X_j - X_i\right), \quad 
		%\omega_{ij} = \frac{1}{2}(\cot\gamma_j + \cot\delta_j) \,, %\quad o = \sum_{j\in N_\triangle(X_i)}A(\tau), 
%\label{eq:N-mesh_1}
%\end{array}	
%\end{equation}	
where $\gamma_j$, $\delta_j$ are the opposite angles to the edge $e_{i,j}$ in the triangles tuple	connected by the edge; %, and $A(\tau)$ denotes the area of the triangle $\tau$.
$D$ is a lumped mass matrix defined as $D=diag\{ |N_\triangle(X_1)|, \cdots ,|N_\triangle(X_n)| \}$ 
which relates to the area/volume around the vertices of the discretized manifold.  

%\begin{lem}
%\label{prop:eigen}

By applying  the discretization $D^{-1}L$ for the LBO on $M$,
and arranging the discretized L$_p$CMs in columns of a matrix $\Psi=\left[ \psi_1, \ldots , \psi_N\right]$, with $\Psi \in \R^{n \times N}$,
the constrained minimization problem (\ref{problem1}) on $M$ reads as follows
\begin{equation}
\label{DP}
\Psi^* \, = \, \arg \min_{\Psi} \frac{1}{\mu} \|\Psi\|_p^p + \mbox{Tr}\left(\Psi^T L \Psi\right) \quad s.t. \quad \Psi^T D \Psi=I,
\end{equation}
where $Tr(\cdot)$ denotes the trace operator, and \tr{ $\|\Psi\|_p^p= \sum_{i,j} d_i |\Psi_{i,j} |^p$, with $d_i$ diagonal elements of the matrix D.}

A discussion on the existence of a minimizer for a constrained variational problem relies on 
conditions on the associated Lagrangian and on the constraints.
In particular, the orthogonality constraints in problem \eqref{DP} are bounded above
by quadratic functions. 
The Lagrangian function of (\ref{DP}) is defined as
\begin{equation}
\label{Lagr}
\mathcal{J}(\Psi,\Lambda) = \frac{1}{\mu} \|\Psi\|_p^p + \mbox{Tr}\left(\Psi^T L \Psi\right) -  \mbox{Tr}\left(\Lambda (\Psi^T D \Psi -I) \right),
\end{equation}
where $\Lambda $ is the matrix of Lagrangian multipliers.  
The function \eqref{Lagr} is proper, lower semi-continuous, bounded from below and coercive. 
If $\Psi$ is a local minimizer of \eqref{DP} then $\Psi$ satisfies the first-order optimality 
conditions 
\begin{equation}
\label{DLagr}
\mathcal{D}_\Psi\mathcal{J}(\Psi,\Lambda) = \frac{1}{\mu} \nu^* + ( 2L - 2 \Lambda D ) \Psi = 0,
\end{equation}
 where 
$\nu^* \;\;{\in}\;\: \partial_\Psi \left[\, \|\Psi\|_p^p \,\right](\Psi^*) $ 
represents the subdifferential (with respect to $\Psi$, calculated at $\Psi^*$), defined in \eqref{diffp},
and we used results from \cite{S85} for trace derivative.

\section{Applying ADMM to the proposed model}
\label{sec:ADMM}

In this section, we illustrate in detail the ADMM-based  iterative algorithm 
used to numerically solve the proposed model (\ref{DP}).
Two different splitting methods for solving problem \eqref{DP} have been proposed in
\cite{OLCO13} and \cite{NVTM}. 
In \cite{OLCO13}  the authors solve the minimization problem by the 
splitting orthogonality constraint (SOC) method introduced in \cite{LO14},
while in \cite{NVTM} an ADMM approach  is introduced that improves the 
empirical convergence performance of the former.
Our approach follows the ADMM strategy, and mainly differs from \cite{NVTM} in the proximal map
sub-problem.   
 
First, we replace the orthogonality constraint  in (\ref{DP}) using an indicator function
$$
\iota(\Psi)=
\left\{
\begin{array}{ll}
0 & \, \mbox{if} \quad \Psi^T D \Psi=I \\
\infty & \mbox{otherwise.}
\end{array}
\right.
$$
Then problem \eqref{DP} can be rewritten as:
\begin{equation}
\label{DPNC}
\Psi^* \, = \, \arg \min_{\Psi}  \frac{1}{\mu} \|\Psi\|_p^p + \mbox{Tr}(\Psi^T L \Psi) + \iota(\Psi).
\end{equation}

We can resort to the variable splitting technique for the orthogonality constraint 
and introduce 
two new auxiliary matrices, $E,S \in R^{n \times N}$, the problem (\ref{DPNC}) is then rewritten as:
\begin{equation}
\label{DPNCSV}
\min_{\Psi,S,E} \frac{1}{\mu} \|S\|_p^p + \mbox{Tr}(E^T L E) + \iota(\Psi) \quad s.t. \quad \Psi=S, \;\, \Psi=E.
\end{equation}
To solve problem (\ref{DPNCSV}), we define the augmented Lagrangian functional
\begin{eqnarray}
\mathcal{L}(\Psi,S,E;U_E,U_S;\mu)
&\;\;{=}\;\;&
\displaystyle{
\frac{1}{\mu} \|S\|_p^p + \mbox{Tr}(E^T L E) +\iota(\Psi)
} \nonumber \\
&&\displaystyle{
{-}\; \langle \, U_S , \Psi - S \, \rangle
\,\;\;{+}\;
\frac{\rho}{2} \, \| \Psi - S \|_F^2
} \nonumber \\
&&\displaystyle{
{-}\; \langle \, U_E , \Psi - E \, \rangle
\;{+}\;
\frac{\rho}{2} \: \| \Psi - E \|_F^2 \;\, ,
} \label{eq:AL}
\end{eqnarray}
where $\rho > 0$ is scalar penalty parameter and $U_S \in \R^{n \times N}$, $U_E \in \R^{n \times N}$
are the matrices of Lagrange multipliers associated with the linear constraints $\Psi = S$ and $\Psi = E$  
in (\ref{DPNCSV}), respectively.

We then consider the following saddle-point problem:
\begin{eqnarray}
\mathrm{Find}&&
\;\;\: (\Psi^*,S^*,E^*;U_S^*,U_E^*)
\;\;{\in}\;\;
\R^{n \times N} {\times}\;\, \R^{n \times N} {\times}\;\, \R^{n \times N} {\times}\;\, \R^{n \times N} {\times}\;\, \R^{n \times N} \nonumber \\
\mathrm{s.t.}&&
\;\;\: \mathcal{L}\,(\Psi^*,S^*,E^*;U_E,U_S;\mu) 
\;\;{\leq}\;\:
\mathcal{L}\,(\Psi^*,S^*,E^*;U_E^*,U_S^*;\mu)
\:\;{\leq}\;\:
\mathcal{L}\,(\Psi,S,E;U_E^*,U_S^*;\mu) \nonumber  \\
&& \;\;\: \forall \: (\Psi,S,E;U_E,U_S)
\;\;\:\!\:\!{\in}\;\;
\R^{n \times N} {\times}\;\, \R^{n \times N} {\times}\;\, \R^{n \times N} {\times}\;\, \R^{n \times N} {\times}\;\, \R^{n \times N}
\: ,
\label{eq:sp}
\end{eqnarray}
with the augmented Lagrangian functional $\mathcal{L}$ defined in (\ref{eq:AL}).

In the following we present the ADMM-based iterative algorithm used to compute a saddle-point 
solution of (\ref{eq:AL})--(\ref{eq:sp}) which 
provides a minimizer of problem (\ref{DP}).

Given the previously computed (or initialized for $k = 0$) matrices
 $S^{(k)}$, $E^{(k)}$, $U_S^{(k)}$ and $U_E^{(k)}$, the $k$-th iteration of the proposed ADMM-based iterative scheme applied to the solution 
of the saddle-point problem (\ref{eq:AL})--(\ref{eq:sp}) reads as follows:
\begin{eqnarray}
&
\Psi^{(k+1)} &
\;{\leftarrow}\;\;\,\,
\mathrm{arg} \: \min_{\Psi \in \R^{n \times N}} \;
\mathcal{L}(\Psi,S^{(k)},E^{(k)};U_S^{(k)},U_E^{(k)})
\label{eq:PM_ADMM_psi} \\
&
S^{(k+1)} &
\;{\leftarrow}\;\;\,\,
\mathrm{arg} \: \min_{S \in \R^{n \times N}} \; 
\mathcal{L}(\Psi^{(k+1)},S,E^{(k)};U_S^{(k)},U_E^{(k)})
\label{eq:PM_ADMM_S} \\
&
E^{(k+1)} &
\;{\leftarrow}\;\;\,\,
\mathrm{arg} \: \min_{E \in \R^{n \times N}} \; 
\mathcal{L}(\Psi^{(k+1)},S^{(k+1)},E;U_S^{(k)},U_E^{(k)})
\label{eq:PM_ADMM_E} \\
&
U_S^{(k+1)} &
\;{\leftarrow}\;\;\,\,
U_S^{(k)} \;{-}\;  \, \rho \, \big( \, \Psi^{(k+1)} \;{-}\; S^{(k+1)} \, \big)
\label{eq:PM_ADMM_US} \\
&
U_E^{(k+1)} &
\;{\leftarrow}\;\;\,\,
U_E^{(k)} \;{-}\;  \, \rho \, \big( \, \Psi^{(k+1)} \;{-}\; E^{(k+1)} \, \big)
\label{eq:PM_ADMM_UE}
\end{eqnarray}

In the following we show in detail how to solve the three minimization sub-problems (\ref{eq:PM_ADMM_psi})--(\ref{eq:PM_ADMM_E}) 
for the primal variables $\Psi$, $S$ and $E$, respectively, \tr{while the ADMM dual variable updates (\ref{eq:PM_ADMM_US})--(\ref{eq:PM_ADMM_UE})
admit closed-form solutions.}
% then we present the overall iterative ADMM-based minimization algorithm.

\subsection{Solution of subproblem (\ref{eq:PM_ADMM_psi}) for $\Psi$}

We observe that the subproblem (\ref{eq:PM_ADMM_psi}) can be rewritten as:
\begin{equation}
\label{one}
\Psi^{(k+1)} \:{\leftarrow}\; \arg \min_\Psi \frac{\rho}{2} \|\Psi- (S + \frac{1}{\rho} U_S) \|_F^2+
\frac{\rho}{2} \|\Psi - (E + \frac{1}{\rho} U_E)  \|_F^2 + \iota(\Psi). 
\end{equation}
%\begin{equation}
%\label{oneb}
%\Psi= \arg \min_\Psi i(\Psi)+\frac{\rho}{2} \|\Psi-\frac{1}{2}(S+E)\|_F^2-2 \mbox{Tr}(S+E) + 
%\frac{1}{2} \|S+E\|_F^2
%\end{equation}
If we omit the constant terms, problem \eqref{one} is equivalent to the following
\begin{equation}
\label{onec}
\Psi^{(k+1)} \:{\leftarrow}\; \arg \min_\Psi \rho \|\Psi-Y\|_F^2  \quad s.t. \quad \Psi^T D \Psi=I
\end{equation}
where $Y=\frac{1}{2} (S + \frac{1}{\rho} U_s + E + \frac{1}{\rho} U_E)$. 

\smallskip

\begin{thm}
\label{th:sol}
The constrained quadratic problem (\ref{onec}), assuming $Y$  has full rank, has the closed-form solution 
\begin{equation}
\label{cfs}
\Psi^{(k+1)} = YV\Sigma^{-1/2}V^T,
\end{equation}
where $V \in \R^{N \times N}$ is a orthogonal matrix and $\Sigma$ is a diagonal matrix satisfying the
SVD factorization $Y^{T}DY  = V\Sigma V^{T}$.
\end{thm}

\smallskip

\begin{proof}
Setting 
\begin{equation}
\label{psi}
\Psi= D^{-\frac{1}{2}}\Phi,
\end{equation}
then the constraint in (\ref{onec}) is equivalent to $\Phi^T \Phi=I$,
and a solution of \eqref{onec} can be obtained by solving:
\begin{equation}
\label{oned}
\min_\Phi \rho \| D^{-\frac{1}{2}}\Phi-Y\|_F^2  \;\; s.t. \;\; \Phi^T  \Phi=I. 
\end{equation}
A closed-form solution of the minimization problem \eqref{oned} 
%\begin{equation}
%\label{minfi}
%\arg \min_\Phi \| D^{-\frac{1}{2}}\Phi-Y\|_F^2 \quad \mbox{s.t.} \quad \Phi^T  \Phi=I
%\end{equation}
can be derived by considering the Lagrangian of the constrained problem \eqref{oned}
\begin{equation}
\label{lagr}
\mathcal{L}(\Phi, \Lambda) = \rho \|D^{-\frac{1}{2}}\Phi-Y\|_F^2+ Tr(\Lambda(\Phi^T\Phi-I))
\end{equation}
where $\Lambda$ is the matrix of Lagrangian multipliers, and its first-order optimality conditions which read as
\begin{equation}
\label{cond}
\left\{
\begin{array}{cl}
\displaystyle{\frac{\partial \mathcal{L}}{\partial \Phi}} = &  2\rho D^{-\frac{1}{2}} (D^{-\frac{1}{2}}\Phi-Y)+\Phi(\Lambda+\Lambda^T)=0\\
\\
\Phi^T\Phi=& I\\
\end{array}
\right.  .
\end{equation}
Multiplying by $D$ the first eq. in \eqref{cond} we obtain:
\begin{equation}
\label{cond1}
\left\{
\begin{array}{rr}
2\rho (\Phi-D^{\frac{1}{2}}Y)+D\Phi(\Lambda+\Lambda^T)= & 0\\
\Phi^T\Phi=& I\\
\end{array}
\right.
\end{equation}
from which it follows that
\begin{equation}
\label{Dy}
D^{\frac{1}{2}}Y=\Phi(I+\hat{D}(\Lambda+\Lambda^T)),
\end{equation}
where $\hat{D}=\frac{1}{2 \rho }D$, and then,
\begin{equation}
\label{fi}
\Phi=  D^{\frac{1}{2}}Y(I+\hat{D}(\Lambda + \Lambda^T))^{-1}
\end{equation}

We set $Z=D^{\frac{1}{2}}Y$, by recalling the second relation of \eqref{cond} and using \eqref{Dy}, it follows
\begin{equation}
Z^TZ=(I+\hat{D}(\Lambda + \Lambda^T))^T(I+\hat{D}(\Lambda+\Lambda^T))
\end{equation}
Since $ Z^T Z \in \R^{N \times N}$, with $N << n$, is symmetric and positive semi-definite, 
following \cite{LO14}, we apply the Singular Value Decomposition (SVD), 
 namely $ Z^T Z = V \Sigma V^T$.
\\
Then $I+\hat{D}(\Lambda+\Lambda^T)=\pm V \Sigma ^{\frac{1}{2}} V^T$ are two square roots of $Z^TZ$. The
principal square root
\begin{equation}
\label{pr}
 (I+\hat{D}(\Lambda+\Lambda^T))=V \Sigma ^{\frac{1}{2}} V^T
\end{equation}
 is the one we desire. If $Z^TZ$ is full rank, then  $V \Sigma ^{\frac{1}{2}} V^T$ is invertible. Thus, relation \eqref{fi} 
can be rewritten as:
\[ %\begin{equation}
\Phi= D^{\frac{1}{2}}Y  V \Sigma ^{-\frac{1}{2}} V^T
\] %\end{equation}
and by \eqref{psi} it follows that
\[ %\begin{equation}
\Psi^{(k+1)}= D^{-\frac{1}{2}}D^{\frac{1}{2}} Y V \Sigma^{-\frac{1}{2}}V^T
\] %\end{equation}
thus \eqref{cfs} holds.
\end{proof}

\noindent{\bf Remark.} \tr{The problem \eqref{oned} is known as orthogonal Procrustes problem. Following \cite{Gow} a solution $\Phi$ of \eqref{oned} reads as 
\begin{equation}
\label{sol_oned}
\Phi=\tilde{U}\tilde{V}^T
\end{equation}
computed by applying the SVD to the matrix $B=(D^{-\frac{1}{2}})^T Y$, thus obtaining $B=\tilde{U}\tilde{\Sigma}\tilde{V}^T $.
Since the SVD computation of an $m \times n$ matrix takes time that is proportional to $O(k m^2 n + k'  n^3)$ with $k$ and $k'$ constants, 
the computational cost for computing the SVD of the $n \times N$ matrix $B$ is $O(n^2N+N^3)$, while 
in the proposed solution, as shown, we computed the SVD of a matrix $Z^TZ$ of dimensions $N \times N$, with a cost of $O(2N^3)$.
Due to the fact that $n >> N$, we conclude that the proposed minimization proved in Theorem \ref{th:sol} is much more computational efficient than the use of the decomposition given in \eqref{sol_oned}.
}

\subsection{Solution of subproblem (\ref{eq:PM_ADMM_S}) for $S$}

Given $\Psi^{(k+1)}, E^{(k)}, U_S^{(k)}$,and $U_E^{(k)}$, and recalling the definition of the augmented Lagrangian functional 
in (\ref{eq:AL}), the minimization sub-problem for $S$ in (\ref{eq:PM_ADMM_S}) can be rewritten as follows:
\begin{equation}
\label{three}
S^{(k+1)} \:{\leftarrow}\; \arg \min_S \frac{1}{\mu} \|S\|_p^p  +\frac{\rho}{2} \|\Psi-(S+\frac{1}{\rho}U_S)\|_F^2
\end{equation}
We can use the Generalized Iterated Shrinkage (GISA) strategy for Non-convex Sparse Coding proposed in \cite{GIS}, 
where the authors extended the popular soft-thresholding operator to $l_p$-norm, or its generalization given
in \cite{LMS15c}.
Rewriting component-wise Eq. (\ref{three}), 
the minimization problem is equivalent to the following $n \times N$ independent scalar problems:
\tr{
\begin{equation}
s_{i,j}^{(k+1)} \:{\leftarrow}\; \mathrm{arg} \min_{s_{i,j} \in \R}
\left\{ \, f(s_{i,j}) = \frac{d_i }{\rho \mu} |s_{i,j}|^p \;{+}\; \frac{1}{2} (s_{i,j} - q_{i,j})^2 \,\right\}
\, , \quad 
\begin{array}{ll}
i = 1, \ldots , n \, , \\
j = 1, \ldots , N \,
\end{array} 
\label{eq:sub_s_i}
\end{equation} }
where $q_{i,j}=\psi_{i,j} - \frac{1}{\rho}(U_S)_{i,j}$.
Following Theorem 1 in \cite{GIS} each of the optimization problems (\ref{eq:sub_s_i}) has
a unique minimum given by
\begin{equation}
\mathrm{prox}_{\frac{d_i}{\rho \mu} \, f}(q_{i,j})
%\mathrm{arg} \min_{s_{i,j} 
\,\;{=}\;
\left\{ \!
\begin{array}{ll}
	0 & 
	\mathrm{if} \;\, | q_{i,j} | \;{\leq}\; \hat{s} \vspace{0.18cm} \\
	\mathrm{sign}(q_{i,j}) \, s^* \;\;\:& 
	\mathrm{if} \;\, | q_{i,j} | \;{>}\; \hat{s} \;\, , 
\end{array}
\right.
\label{eq:prox_phi_sol}
\end{equation}
where the thresholding value is 
\tr{
\[
\hat{s}= \left( \frac{2d_i}{\rho \mu}(1-p) \right)^{1/(2-p)}+\frac{d_i}{\rho \mu}p \left( \frac{2d_i}{\rho \mu}(1-p) \right)^{(p-1)/(2-p)}
\] 
}
and $\,s^*$ is the unique solution of the following nonlinear equation:
\tr{
\begin{equation}
s_{i,j} - q_{i,j} +p\frac{d_i}{\rho \mu} (s_{i,j})^{p-1} \, = \, 0 \;,
\end{equation}
}
that can be easily solved by a few iterations of an iterative zero-finding algorithm.

\subsection{Solution of subproblem (\ref{eq:PM_ADMM_E}) for $E$}
Given $\Psi^{(k+1)}, S^{(k+1)}, U_S^{(k)}$, and $U_E^{(k)}$, the minimization problem of the augmented Lagrangian functional 
in (\ref{eq:AL}) with respect to $E$ in (\ref{eq:PM_ADMM_E}) can be rewritten as follows:

\begin{equation}
\label{two}
E^{(k+1)} \:{\leftarrow}\; \arg \min_E Tr(E^T L E) +\frac{\rho}{2} \|\Psi-(E+\frac{1}{\rho}U_E)\|_F^2
\end{equation}
To solve the minimization problem \eqref{two}, we consider the optimality conditions, namely:
$$2 L E+ \rho (\Psi-(E+\frac{1}{\rho}U_E))=0$$
which reduce to the solution of $N$ linear systems for $E$  in the following form
\begin{equation}
\label{SolE}
(\rho I - 2 L)E= \rho(\Psi-\frac{1}{\rho} U_E).
\end{equation}

\section{Basics on partitioning}
\label{sec:basic}

The $N$ orthogonal L$_p$CMs have the potential to be localized in the $N$ main key features of the shape.
This can be naturally exploited to subdivide the shape into a collection of salient parts.

Shape partitioning enables the decomposition of arbitrary topology objects 
into smaller and more manageable pieces called partitions.
In particular we are interested in Manifold Partitioning, since the boundaries of tangible physical objects 
can be mathematically defined by two-dimensional manifolds embedded into three-dimensional Euclidean space.

Let us introduce the following formulation of the shape partitioning problem.

\bigskip

\begin{defn}[Manifold Partitioning]\label{def:MP}
Given a compact 2-manifold $\mathcal{M}$, find the partition into $N$ sub-manifolds defined by the pairs of topological spaces 
$\{ (U_k, \partial U_k ) \}_{k=1}^N$, with boundary $\partial U_k$,  
such that all of the following conditions hold:
\begin{itemize}
	\item[P1)] $U_k, k=1,\ldots,N,$ is a non-empty connected sub-manifold; %does not contain  empty set.
	\item[P2)] $\displaystyle{ \bigcup_{k=1}^N U_k = {\mathcal M}}$ 
	%The union of the sets in P is equal to X. (The sets in P are said to cover X.)
  \item[P3)] The intersection of any two distinct sub-manifolds $U_i, U_j$ in $\mathcal{M}$ is equal to a simple curve:
	 \[{\displaystyle U_i \cap U_j = \partial U_i \cap \partial U_j =}  \mbox{1-manifold}.\]
%	 \item[P4)] $U_k$ is a connected sub-manifold, $k=1,\ldots,N$.\\
%	(We say the elements of P are pairwise disjoint
\end{itemize}
\end{defn}

\bigskip

The sub-manifolds $\{U_k\}_{k=1}^N$ are said to cover $\mathcal{M}$ and provide the so-called \textit{segmentation}, or partitioning, 
 of the object represented by ${\mathcal M}$.

Many shape processing applications rely on a more stringent characterization of partitioning 
which requires a global parametrization of the manifold.
However, smooth global parameterization does not always exist or is easy to find.
Only the simplest 2-manifolds indeed can be adequately parameterized.
In general, a topology decomposition of the manifold is required to describe it as a
a collection of parameterized surfaces (charts).

%A parameterized surface can be described by the map                                                
%$S: D \subset \R^2 \rightarrow \R^3$  is a parameterization domain.
%the map  $S(u,v):=(x(u,v),y(u,v),z(u,v))$                                  
%is a global parameterization (embedding) of   the manifold.
%
%
%Similarly, a differentiable manifold can be described using mathematical maps, called coordinate charts, collected in a mathematical atlas. 
%It is not generally possible to describe a manifold with just one chart, because the global structure of the manifold is different from the simple structure of the charts.
%An atlas is not unique as all manifolds can be covered multiple ways using different combinations of charts.
We \tr{briefly} review some useful definitions.

A chart for a 2-manifold  ${\mathcal{M}}$ is a homeomorphism ${\displaystyle \varphi }$  from a subset $U$ of 
 ${\mathcal{M}}$ to a subset of the two-dimensional Euclidean space. 
The chart is traditionally recorded as the ordered pair ${\displaystyle (U,\varphi )}$.
A collection ${\displaystyle \{(U_{k},\varphi_{k})\}}$ of charts on ${\mathcal{M}}$
 such that ${\displaystyle \bigcup U_{k}=\mathcal{M}}$ forms an atlas for ${\mathcal{M}}$. 

When a manifold is constructed from multiple overlapping charts, the regions where they overlap carry information essential for understanding the global structure. In this context, as specified by $\mathit{P3)}$ in Definition \ref{def:MP}, the 
overlap is reduced to boundary curves shared by two adjacent patches.

A \textit{patch-based} partitioning can be then defined as follows.

\bigskip

\begin{defn}[Patch-Based Manifold Partitioning]\label{def:PB_MP}
Given a compact 2-manifold $\mathcal{M}$, find the partition into $N$ sub-manifolds $\{ (U_k, \partial U_k ) \}_{k=1}^N$
such that conditions P1) - P3) hold, together with the following 
\begin{itemize}
	\item[P4)] $U_k$ is a genus-0 sub-manifold that defines a chart.
	\item[P5)] $U_k$ has at most two boundaries.
\end{itemize}
\end{defn}
%%%

\bigskip

Given a chart decomposition of a mesh, each chart can be parameterized on a planar domain
(e.g., a circle or a rectangle) using different methods, whose selection depends on its genus and
number of boundary components. More precisely, a disk-like charts (i.e., genus-0 patches with one boundary component) are parameterized using the barycentric coordinates method \cite{Floater2005};
while a genus-0 chart with more than one boundary component, or more generally charts with an arbitrary
genus, are converted to disk-like regions by cutting them along cut-graphs and then embedded on the
plane using the barycentric coordinates method \cite{PATANE2007},\cite{STEINER2002}.

For approximation purposes and in order to reduce the parameterization distortion, it is preferable
to work with disk-like patches.

In \cite{PSF2004} a topology-based decomposition of the shape is computed and used to segment the shape into
primitives, which define a chart decomposition of the mesh. The charts considered in \cite{PSF2004} are all genus-0 
but can present more than one boundary components. 
\tr{In contrast}, in this work we restrict the chart $U_{k}$ to be 
a disk-like patch bounded by one or two closed curves. 
The latter requires a simple cut between the two boundaries to avoid internal holes in the planar parameterization.

Once the proposed patch-based manifold partitioning is built, we can associate a parameterization $\varphi_{k}$
to each sub-manifold $U_{k}$. However, we omit the construction of a parameterization, 
as discussing these details goes beyond the scope of this paper.

\section{The Partitioning  Algorithm}
\label{sec:alg}

In Section \ref{sec:ADMM} we described an optimization method to compute a basis of $N$ functions $L_p$CMs 
induced by the LBO of a manifold $\mathcal{M}$ represented by a mesh $M$ with $n$ vertices.
Each $L_p$CM has compact support: it is non-zero only in a confined region of the domain,
and the size of the compact support can be controlled by $\mu$ and $p$.
 
We propose a numerical algorithm to partition a mesh which iteratively increases 
the support  of $N$ functions $L_p$CMs, 
with $N <<n$, until their supports cover the entire mesh without overlapping.
The set of vertices in the support of $\psi_i$ defines a sub-mesh.
A partitioning of $M$ is defined as the \tr{union} of the $N$ sub-meshes $\psi_i, i=1,\ldots,N$.

The algorithm consists of three main steps illustrated in \textbf{Algorithm 1},
which takes as input the initial mesh $M$, the number of partitions $N$ or the initial $\mu$ value,
and returns a set of sub-meshes $S$. 
As concerning mesh segmentations, given in Def. 1, the set $S$ is directly the output of Step 2,
while for patch-based partitioning a further step (Step 3) is required to suitably refine the partition $S$ 
according to Def. 2.
%to satisfy $P4)$ in (\ref{def:PB_MP}).     

\bigskip

%The presented ideas are summarized in the following pseudo-code \textbf{Algorithm 1}.

%%%%
\begin{algorithm}
	\label{alg:1}
	\noindent\caption{\bf Mesh Partitioning \vspace{0.05cm} }
	\vspace{0.2cm}
	{\renewcommand{\arraystretch}{1.0}
		\renewcommand{\tabcolsep}{0.0cm}
		\vspace{-0.08cm}\\
		\begin{tabular}{ll}
			\textbf{Input}:       & mesh $M \,$, $\mu$ or $N$\vspace{0.04cm} \\
			\textbf{Output}: & patch set $\,S=\{S_k\}_{k=1}^N$  \vspace{0.2cm} \\
			\textbf{Parameters}: $\;\;$& tolerance $\epsilon=0.01$  \vspace{0.0cm} \\
			\vspace{0.0cm}
		\end{tabular}
	}\\
	
	\vspace{0.1cm}
	
	{\renewcommand{\arraystretch}{1.0}
		\renewcommand{\tabcolsep}{0.0cm}
		\begin{tabular}{lll}
			\hline
			\multicolumn{3}{l}{STEP 1: \hspace{0.05cm}Compute ${\Psi}\in \R^{n\times N}$} \vspace{0.05cm}\\
			\hline
			& &  \\
			%\vspace*{0.1cm}
			\hspace*{0.05cm}& 
			\begin{tabular}{ll}		
				STEP 1a (given $\mu$): &  \\
				$\quad\cdot$ set $uncovered$  = true, $N=1$ & \\
				\multicolumn{2}{l}{$\quad$\textbf{while} $\;$ ($uncovered$)  \textbf{do}:} \vspace{0.1cm}\\
				$\quad\quad\cdot$ $N \leftarrow N+1 $ \vspace{0.1cm} &\\
				$\quad\quad\cdot$ Compute $\{\psi_i\}_{i=1}^{N} \;$ by solving  (\ref{DP}) \vspace{0.1cm} &\\
				$\quad\quad\cdot$ set $uncovered  \, = \, (\exists X_j : \psi_i(X_j) \, = \, 0 \; \forall i)$ & \\
				\multicolumn{2}{l}{\textbf{$\quad$end$\;$while}} \vspace{0.09cm} \\
			\end{tabular}				& 
			\begin{tabular}{ll}		
				STEP 1b (given $N$): & \\
				$\quad\cdot$ set $uncovered$  = true, $\mu=2$ &\\
				\multicolumn{2}{l}{$\quad$\textbf{while} $\;$ ($uncovered$)  \textbf{do}:} \vspace{0.1cm}\\
				$\quad\quad\cdot$ update $\mu \leftarrow 4\mu $ \vspace{0.1cm} &\\
				$\quad\quad\cdot$ Compute $\{\psi_i\}_{i=1}^{N} \;$ by solving  (\ref{DP}) \vspace{0.1cm} &\\
				$\quad\quad\cdot$ set $uncovered  \, = \, (\exists X_j : \psi_i(X_j) \, = \, 0 \; \forall i)$ & \\
				\multicolumn{2}{l}{\textbf{$\quad$end$\;$while}} \vspace{0.09cm} \\
			\end{tabular} \\
			\hline
			\vspace*{-0.3cm} & &  \\
			\multicolumn{3}{l}{STEP 2: \hspace{0.05cm}Region Growing} \vspace{0.05cm}\\
			\hline
			& &  \\
			%%
			%& \multicolumn{2}{l}{\textbf{for} $\;$ \textit{k = 1, $\, \ldots \, N$ $\:$} \textbf{do}:} \vspace{0.1cm}\\
			%
			%& $\quad\cdot$ set seeds $s_k$ according to (\ref{eq:seeds}), $\quad S_k = s_k$ & \\
			%& $\quad\cdot$ set initial buffer $b_k \in B\,,\quad b_k = N_\triangle(s_k)$ & \\
			%& \multicolumn{2}{l}{\textbf{end$\;$for}} \vspace{0.09cm} \\
			%
			%& \multicolumn{2}{l}{\textbf{repeat} $\;$  }\vspace{0.1cm}\\
			& \multicolumn{2}{l}{$\quad$\textbf{for} $\;$ \textit{k = 1, $\, \ldots \, N$ $\:$} \textbf{do}:} \vspace{0.1cm}\\
			& $\quad\cdot$ set seeds $s_k$ according to (\ref{eq:seeds}),  & \\
			& $\quad\cdot$ set initial buffer  $b_k \leftarrow N_\triangle(s_k)$ & \\
			& \multicolumn{2}{l}{$\quad\quad$\textbf{while} $\;$ $b_k \neq \{ \oslash \}$  \textbf{do}:} \vspace{0.1cm}\\
			& $\quad\quad$ \textbf{if} $\;$ $(||{\psi}_k(\tau)| - \max\limits_{i=1,..,N}|{\psi}_i(\tau)|| \leq \epsilon )$ \vspace{0.05cm} &\\
			& $\quad\quad\quad\cdot$ add $\:\tau$ in $\;S_k$ \vspace{0.1cm} &\\
			& $\quad\quad\quad\cdot$ update $\,b_k$ by inserting $\,N_\triangle(\tau)\;$ \vspace{0.1cm} &\\
			& \multicolumn{2}{l}{\textbf{$\quad\quad$end$\;$if}} \vspace{0.09cm} \\
			& $\quad\quad\quad\cdot$ update $\,b_k$ by removing $\tau$ \vspace{0.1cm} &\\
			& \multicolumn{2}{l}{\textbf{$\quad\quad$end$\;$while}} \vspace{0.09cm} \\
			& \multicolumn{2}{l}{\textbf{$\quad$end$\;$for}} \vspace{0.09cm} \\
			%& \multicolumn{2}{l}{\textbf{until $| B | > 0$}}    \vspace{0.09cm} \\
			%& \multicolumn{2}{l}{$u^* = u^{(k)}$}\\
			\hline
			\vspace*{-0.3cm} & &  \\
			\multicolumn{3}{l}{STEP 3: \hspace{0.05cm} Refinement for Patch-Based Manifold Partitioning }\\
		\end{tabular}
	}
\end{algorithm}
%%%%

In Step 1 an iterative process is applied to generate
 a basis $\{\psi_i\}_{i=1}^{N}$ by solving (\ref{DP}) with the ADMM procedure described in Sec.\ref{sec:ADMM}.
This task can be realized following two different approaches, named Step 1a and Step 1b,
that terminate when all the vertices $V$ are in the support of at least one L$_p$CM.

In Step 1a the parameter $\mu$ in (\ref{DP}) is assigned. Starting from the 
construction of a small set of functions $L_p$CMs, the space dimension is 
enlarged at each iteration by adding a new function $\psi_{i}$
until any vertex of $M$ is covered by at least one function in $\Psi$.
Starting from a small dimension space, Step 1a ends up with 
a space of dimension $N$ spanned by the $L_p$CMs. 
In this approach the final number of partitions $N$ is unpredictable in advance.
 
Alternatively, Step 1b overcomes the problem to identify an a priori value for $\mu$
and requires instead a fixed number for $N$. At each iteration, $N$ basis functions $L_p$CM are built by solving (\ref{DP}) with a given 
$\mu$. If there exists a vertex of $M$ not covered by any function in $\{\psi_i\}_{i=1}^{N}$, then 
$\mu$ is increased, thus causing an enlargement of the function supports. The 
solution of (\ref{DP}) is then re-iterated with the new value for $\mu$.

Once the $N$ functions discretized in $\Psi\in\R^{n\times N}$ are determined by either Step 1a or Step 1b, 
the whole set of vertices $V$ is covered but many regions can be over-covered by more L$_p$CMs. 
Mesh partitioning satisfying Def. 1 is then carried out in Step 2. At this aim, 
$N$ initial seeds $(s_1,\ldots,s_N)$ are selected as the $L_p$CM extrema,
as follows 
\begin{equation}
\label{eq:seeds}
s_k = \arg \max\limits_{i=1,..,n}\left|{\psi}_k (X_i) \right|\, \quad k=1,\ldots,N.
\end{equation}

%Once the $d$ Compressed Modes $\Psi\in\R^{n\times d}$ are obtained, 
%such that the whole set of vertices $V$ is covered, 
%we can use them as a descriptor for the patching algorithm.
%The algorithm is based on the region growing triangles 
%with $d$ initial seeds $s_k$ placed in the CM extrema
%\begin{equation}
%\label{eq:seeds}
%s_k = \arg \max\limits_{i=1,..,n}\left|\bar{\psi}_i^k\right|\,,
%\end{equation}
%where the function set $\bar{\Psi}\in\R^{n_T\times d}$ defined over the mesh triangles is obtained
%easily by interpolating $\Psi$ at triangles barycentre.
%Seeds selected in this way are distinct thanks to the local support of the compressed modes
%and the mechanism how they were constructed.

Then a region growing strategy is applied which consists of a buffer of adjacent neighbors of a given element set, and a loop
in which the buffer and the element set are updated according to some decision rule.
%Our goal is to split the input triangular mesh into $d$ patches defined by the dimension 
%of the L$_p$CM model output.
Starting from the initial buffer $b_k = N_{\triangle}(s_k), k=1,\ldots,N,$
we examine each triangle $\tau\in b_k$ to decide if it will be added to $S_k$
which is the set of triangles associated with the function ${\psi}_k$. 
We denote by ${\psi}_k(\tau)$ the value obtained interpolating ${\psi}_k$ on its vertices.

There are two cases that may occur when an  unassigned triangle  $\tau$ is considered: 
\begin{itemize}
	\item In the first case, $\tau$ is covered by one support, 
	 i.e. ${\psi}_k(\tau) \neq 0\,$ and ${\psi}_i(\tau) = 0 \, \; \forall i\ne k$.
We remove $\tau$ from $b_k$ and assign it to $S_k$.
Then the buffer $b_k$ is updated by adding the $\tau$'s neighbors $N_\triangle(\tau)$.
	\item The second case occurs when the supports of at least two basis functions overlap, i.e. 
				${\psi}_k(\tau) \neq 0\,$  and $\exists\, i\ne k \, : \, {\psi}_i(\tau) \ne 0\,$.
This case locates over the bands of overlapped supporting functions.
%A natural selection rule would be to assign $\tau$ to set $S_*$ associated with 
%a supporting function highest in magnitude $* = \arg\max\limits_{i=1,..,d}|\bar{\psi}_i(\tau)|$.
%In practice we experienced that such extrema decision rule is not robust enough in a narrow band 
%where the supporting functions are decreasing towards zero.
%Due to occasional oscillations $\tau$ was misplaced what caused the algorithm to fail.
%Therefore, we decided rather to 
If the difference from the extrema is under a certain threshold $\epsilon$,
which reads as
	\begin{equation}
\left| |\psi_k(\tau)| - \max_{i=1, \ldots , N}|{\psi}_i(\tau)| \right| \leq \epsilon,
	\label{eq:max}
	\end{equation}
then $\tau$ is added to $S_k$ and $b_k$ is updated accordingly, as in the previous case.
Otherwise, 
%	\item The third case is complementary to the latter, what means 
$\tau$ will be assigned to  a different set and 
% does not belong in $S_k$ and it will be included in a different set.
the only action taken in this case will be to remove $\tau$ from $b_k$.
\end{itemize}
We notice that condition (\ref{eq:max}) is trivially satisfied in the first case. 

\smallskip

\begin{figure*}[ht]
	\centering
	\includegraphics[width=5.1cm]{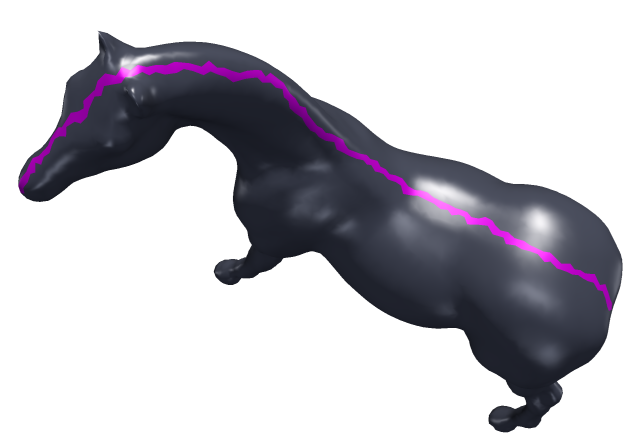}	\hfill
	\includegraphics[width=7.1cm]{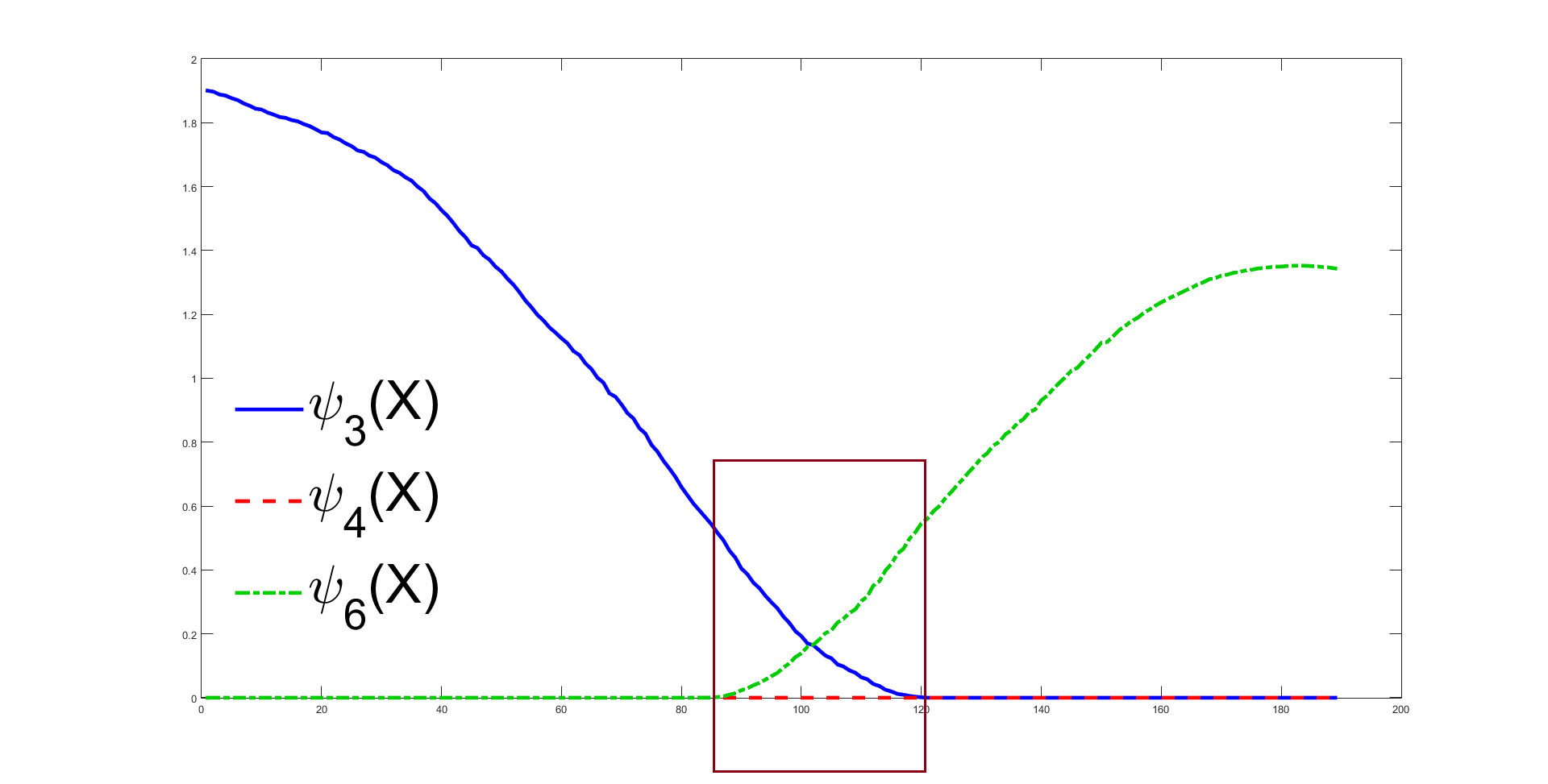}	\\	%\hspace{-0.5cm}
	
	\caption{ L$_p$CMs plotted over the magenta curve on the \texttt{horse} mesh.
		Over the curve only $\psi_3$ and $\psi_6$ are non-zero, and the red box shows the 
		band where the supports overlap.}
		%Function $\psi_4$ that localizes rear-left leg is plotted as justification of the local support 
		%since over the highlighted line takes constant zero value.}
	\label{fig:CMline}
\end{figure*}

A better understanding of condition (\ref{eq:max}) is provided in Fig.\ref{fig:CMline}
The region growing step has been applied to partition the \texttt{horse} mesh into $N=6$ parts;
 the partitioning results of Step 2 are illustrated in Fig.\ref{fig:patches}.
	
Along the magenta colored curve depicted on the mesh (Fig.\ref{fig:CMline}, left), from the horse's head to its bottom, 
we plot the values of the L$_p$CMs (Fig.\ref{fig:CMline}, right).
Only the two functions $\psi_3$ and $\psi_6$ of $\Psi$ are non-zero.
For the sake of clarity we plot also  $\psi_4$, which localizes rear-left leg of the \texttt{horse} mesh
	and over the line evaluates zero.
The red box locates the band of overlapping.
When the functions values, e.g. $\psi_3$ and $\psi_6$, are too close (below $\epsilon$), 
	even in case of some minor numerical perturbation, a corresponding set of successive triangles 
	may tend to \tr{over-leap} in the cluster assignment.
However, the condition (\ref{eq:max}) satisfyingly overcomes this practice issue.

Step 2 ends when the buffers are empty, i.e. all the triangles of M have been assigned to $S$.

An a posteriori procedure approximates the boundaries of the 
sub-meshes $\{S_k\}_{k=1}^N$ by smooth spline curves.

Step 3 of \textbf{Algorithm 1} is applied to finalize the Patch-Based Manifold Partitioning
following Definition 2.
The refinement is required only for a few patches $\{ S_k\}_{k=1}^{\bar{N}}$, $\,  \bar{N}<N$,
with genus greater than
zero, and for genus-0 patches with more than two closed-loop boundaries.
%The genus of a patch is computed by applying the Euler’s formula for mesh with boundaries. 
%A patch of genus higher than 0 is easily defined by the number of its closed-loop boundaries.
%If a patch has more than two closed boundaries, there would be needed at least two cuts along the patch
%to create two separate parts.
The refinement is an adaptive process that consists in the re-iteration of Step 1 and Step 2 for every patch 
$S_k$ that needs to be further subdivided, by imposing the initial number of partitions $N=2$.

\section{Experimental Results}
\label{sec:ne}

In this section we describe the experimental results which demonstrate the performance of \textbf{Algorithm 1}. 
In particular, we first evaluate the performance of Step 1 for the computation of the $L_p$ Compressed  
Modes $\Psi\in\R^{n\times N}$, then we illustrate the results of Step 2 and Step 3 for part-/patch-based partitioning, respectively. 

Experimental tests were performed on Intel\textregistered Core\texttrademark i7-4720HQ 
Quad-Core 2.6 GHz machine, with 12 GB/RAM and Nvidia GeForce GTX 860M graphics card 
in a Windows OS.
The code is written in \textsc{Matlab}, and executed without any additional machine support, e.g. parallelization and GPU-based computations.

We tested the proposed method on a set of meshes downloaded 
from the data repository website \texttt{http://segeval.cs.princeton.edu}, \cite{Chen}. 
The dataset represents  geometric models with different characteristics in terms of details, level of refinement, and present a medium dense vertex distribution, in particular the number of vertices
	\tr{and triangles}
	of the meshes visualized in the examples \tr{are reported in the second and third column of Table \ref{tab:Data}}.

The figures reported in this section were produced by the software ParaView, and its VTK reader.
In the examples illustrated we applied a post-process smoothing to the boundaries between 
the segmented parts $\{S_k\}_{k=1}^N$ by projecting the boundary vertices onto the cubic spline obtained 
by least-squares approximation.

\subsection{STEP 1: Computing the L$_p$CMs }

The two strategies Step 1a and Step 1b in \textbf{Algorithm 1}, described in section \ref{sec:alg},
	generate the basis functions $\Psi$.

In all the experiments we used a randomized matrix as initial iterate $\Psi^{(0)}$ for the ADMM computation of (\ref{DP}), 
and we terminated the ADMM iterations as soon as the relative 
change between two successive iterates satisfies 
\begin{equation}
\, err_{\Psi} = \frac{\| \Psi^{(k)} - \Psi^{(k-1)} \|_{F}}{\| \Psi^{(k-1)}\|_{F}} \;{<}\; 10^{-3}.
\label{errpsi}
\end{equation}
%The algorithm is initialized with a random $\Psi$ matrix. 
As already observed in \cite{NVTM}, where the L$_1$ penalty term is used, 
different runs converge to the same set of basis
functions, although their ordering might be different. 
\tr{In our experiments the $p$ values were tested in the range $[0.5,0.8]$.
However, since small $p$ values affect mainly the efficiency, we decided to set the sparsity parameter $p=0.8$ for all the examples reported. }

\begin{figure*}[ht]
	\flushleft
	\includegraphics[width=2.05cm]{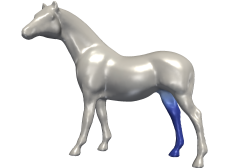} \hspace{-0.5cm}
	\includegraphics[width=2.05cm]{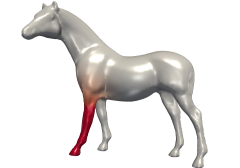} \hfill
	\includegraphics[width=2.05cm]{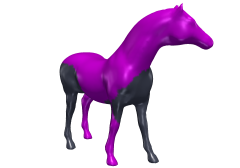} \\
	\includegraphics[width=2.05cm]{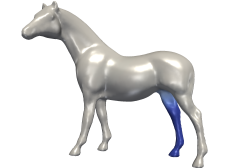}	\hspace{-0.5cm}
	\includegraphics[width=2.05cm]{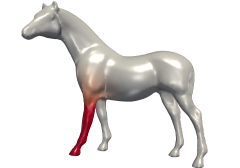} \hspace{-0.5cm}
	\includegraphics[width=2.05cm]{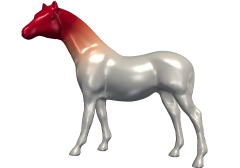} \hfill
	\includegraphics[width=2.05cm]{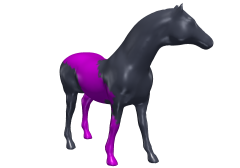} \\
	\includegraphics[width=2.05cm]{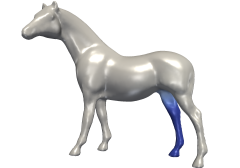} \hspace{-0.5cm}
	\includegraphics[width=2.05cm]{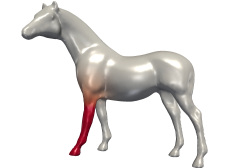} \hspace{-0.5cm}
	\includegraphics[width=2.05cm]{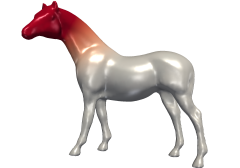} \hspace{-0.5cm}
	\includegraphics[width=2.05cm]{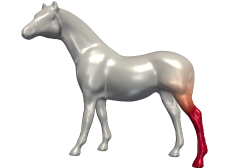} \hfill
	\includegraphics[width=2.05cm]{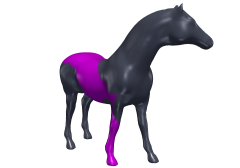} \\
	\includegraphics[width=2.05cm]{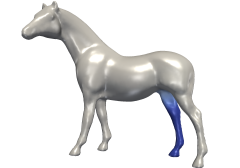} \hspace{-0.5cm}
	\includegraphics[width=2.05cm]{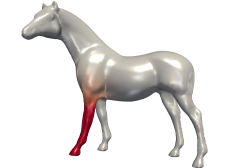} \hspace{-0.5cm}
	\includegraphics[width=2.05cm]{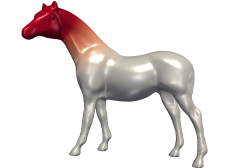} \hspace{-0.5cm}
	\includegraphics[width=2.05cm]{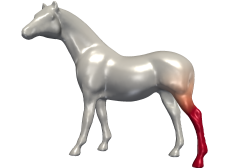} \hspace{-0.5cm}
	\includegraphics[width=2.05cm]{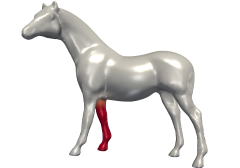} \hfill
	\includegraphics[width=2.05cm]{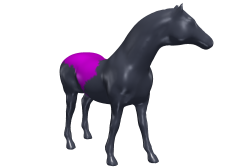} \\
	\includegraphics[width=2.05cm]{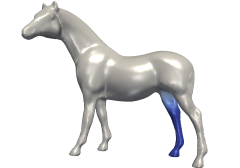} \hspace{-0.5cm}
	\includegraphics[width=2.05cm]{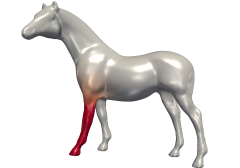} \hspace{-0.5cm}
	\includegraphics[width=2.05cm]{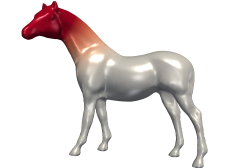} \hspace{-0.5cm}
	\includegraphics[width=2.05cm]{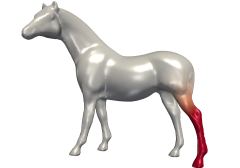} \hspace{-0.5cm}
	\includegraphics[width=2.05cm]{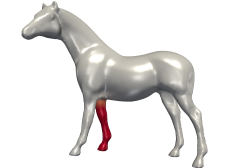} \hspace{-0.5cm}
	\includegraphics[width=2.05cm]{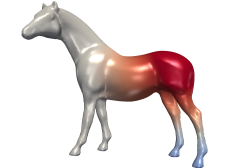} \hfill
	\includegraphics[width=2.05cm]{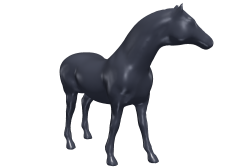} 
	
	\caption{ L$_p$CMs generated at each iteration of Step 1a in Algorithm 1, on the \texttt{horse} mesh.}
	\label{fig:V1_scheme}
\end{figure*}

Figure \ref{fig:V1_scheme} illustrates how Step 1a works when the parameter value $\mu$ is fixed, $\mu=300$.
At the first iteration, only two initial quasi-eigenfunctions are computed with the given $\mu$. 
The control of the local support volume resulted in localizing \tr{two legs} of the \texttt{horse} 
	mesh, leaving the rest uncovered (highlighted in magenta at the end of the first row).
In the second iteration (second row), the space dimension is enlarged ($N=3$), resulting in optimization
	of $\Psi^{n\times 3}$.
The support of the third function $\psi_3$ shrinks the uncovered area under the head and neck, leaving just two legs and 
part of the \texttt{horse}'s body uncovered.
%In general, due to the random-values initialization, the first two functions can be localized
	%at different parts of the mesh for the same given $\mu$.
%Therefore at each run, can cover the mesh parts in different order.
%Notice that in at this point the third supporting function is non-zero over both the rear legs.
%Even if the values of $\psi_3$ at this iteration seem disconnected, we can see on the right that
%	the whole bottom part of the mesh is covered by this supporting function, leaving just one leg
%	and the horse's belly uncovered.
The algorithm terminates after five iterations, enlarging the space up to six functions
	$\psi_1,\dots,\psi_6$ and leaving no more vertices of $M$ uncovered.
The result of the last iteration is depicted in the bottom row of Figure \ref{fig:V1_scheme}.
Notice that over iterations, the corresponding functions describing the same parts of the mesh
	retain their order in the set $\Psi$.
%This, in general, does not hold.
Due to the randomized initialization, the order in general changes for different runs.

\begin{figure*}[ht]
	\centering
	\includegraphics[width=3cm]{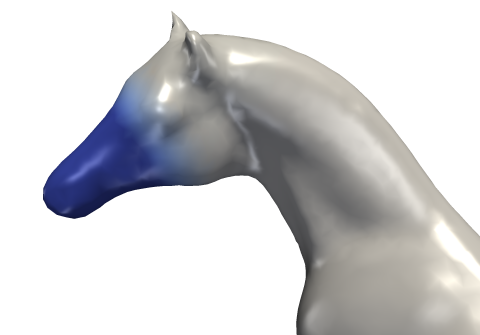}	\hfill
	\includegraphics[width=3cm]{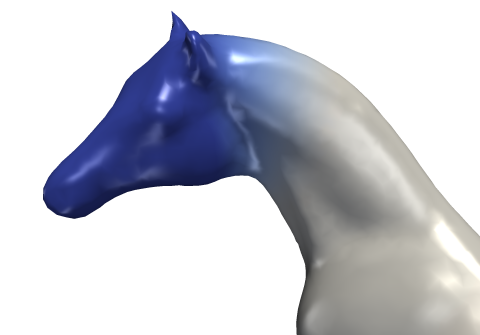}	\hfill
	\includegraphics[width=3cm]{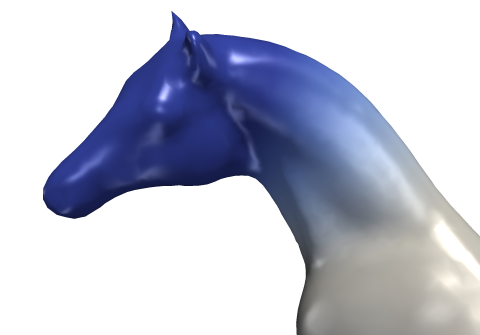}	\hfill
	\includegraphics[width=3cm]{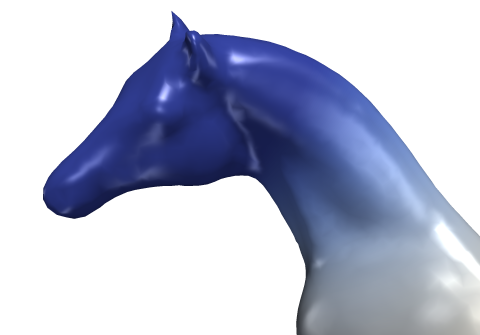}	\\
	\includegraphics[width=3cm]{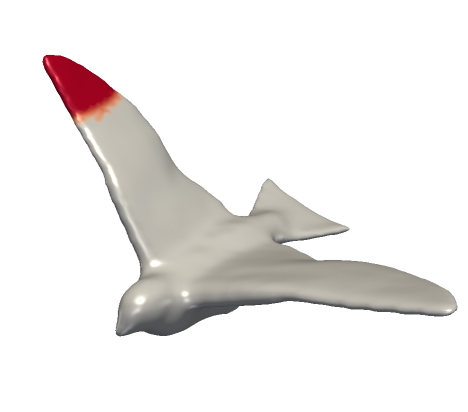}	\hfill
	\includegraphics[width=3cm]{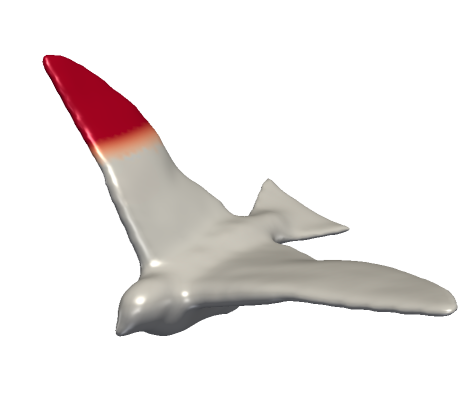}	\hfill
	\includegraphics[width=3cm]{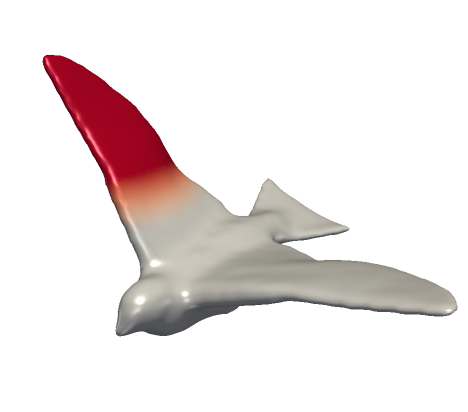}	\hfill
	\includegraphics[width=3cm]{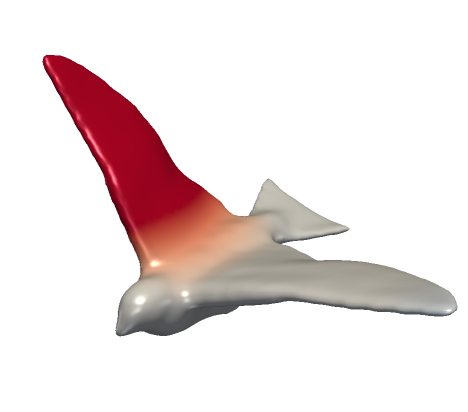}	\\
	
	\caption{L$_p$CM $\psi_5$ of $\Psi$ \tr{on the \texttt{horse} mesh (top), and $\psi_1$ of $\Psi$ on the \texttt{bird} mesh (bottom), both} generated for 4 iterations of Step 1b in Algorithm 1.
	From left to right: enlarging of the support obtained by increasing the parameter $\mu$ for a fixed number of basis functions.}
	\label{fig:V2_scheme}
\end{figure*}

Step $1b$ iteratively recomputes a given number $N$ of basis functions increasing the value of the $\mu$
parameter, thus enlarging the local support at each iteration, until all the vertices of $M$ are covered by at least one function $\psi$.
By the way of illustration, in Figure \ref{fig:V2_scheme} we show the enlargement of the support of $\psi_5 \in \Psi$ for \tr{\texttt{horse} mesh and $\psi_1 \in \Psi$ for \texttt{bird} mesh, for}  
increasing values of $\mu$ and a fixed basis dimension $N=5$ \tr{and $N=4$ respectively}.
%We run \texttt{STEP 1} of the \textbf{Algorithm 1} asking for a set of functions of dimension 5.
The initial $\mu=8$ is increased by a factor 4 at each iteration.
From left to right, the results are shown for $\mu=8$, $\mu=32$, $\mu=128$ and $\mu=512$\,.
%

%In figure \ref{fig:CM_mu_dep} we demonstrate the effect of the parameter $\mu$ on the volume
%	of the support of the L$_p$CMs. 
%Applying the ADMM-based algorithm to problem (\ref{DPNC}), enlarging the space dimension 
%	for a fixed $\mu$ until the basis covered the whole mesh $M$.
%Upper row shows basis of 7 functions computed for $\mu=0.008$, 
%	the middle row covers the mesh with 5 basis functions for $\mu=0.003$
%	and in the bottom row the value $\mu=0.001$ used produced just 3 basis functions,
%	however, we cannot describe them as basis functions with local support anymore.
%The uncoloured parts of the mesh are located outside the support of each basis function, i.e.
%	the value $\psi_i(X_j) = 0$.

\begin{figure*}[ht]
	\centering
	\includegraphics[width=2.5cm]{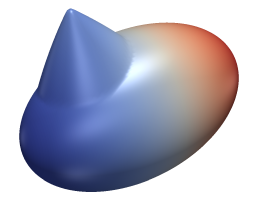}	\hfill
	\includegraphics[width=2.5cm]{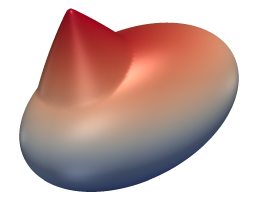}	\hfill
	\includegraphics[width=2.5cm]{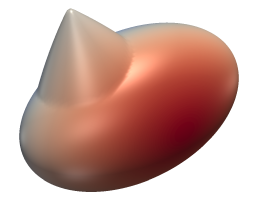}	\hfill
	\includegraphics[width=2.5cm]{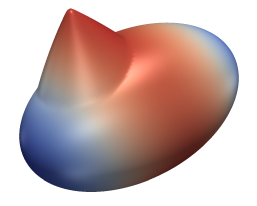}	\hfill
	\includegraphics[width=2.5cm]{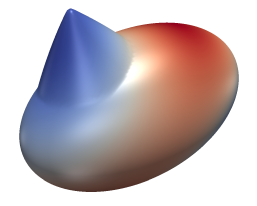}	\\
	\includegraphics[width=2.5cm]{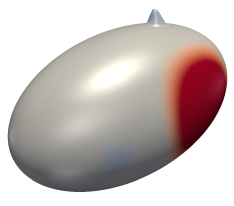}	\hfill
	\includegraphics[width=2.5cm]{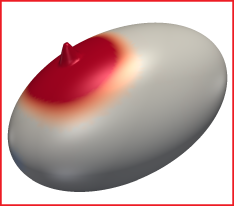}	\hfill
	\includegraphics[width=2.5cm]{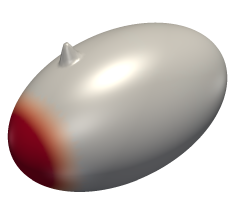}	\hfill
	\includegraphics[width=2.5cm]{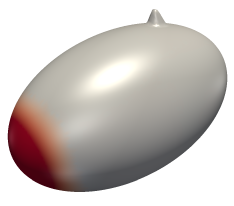}	\hfill
	\includegraphics[width=2.5cm]{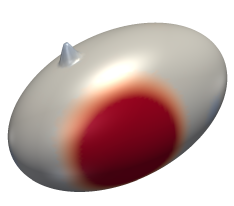}	\\
	\includegraphics[width=2.5cm]{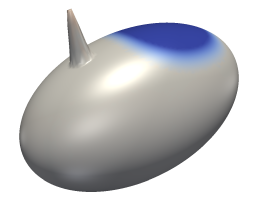}	\hfill
	\includegraphics[width=2.5cm]{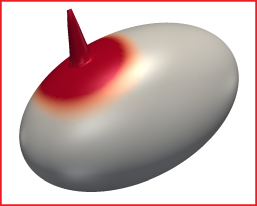}	\hfill
	\includegraphics[width=2.5cm]{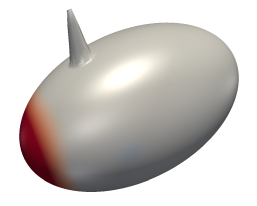}	\hfill
	\includegraphics[width=2.5cm]{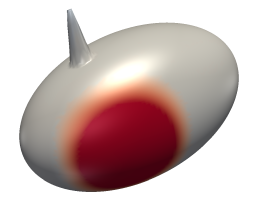}	\hfill
	\includegraphics[width=2.5cm]{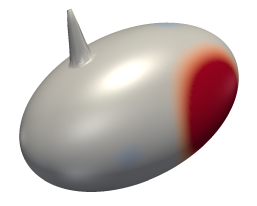}	\\
	\includegraphics[width=2.5cm]{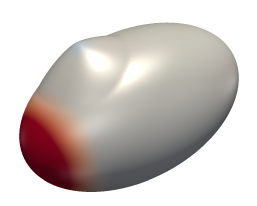}	\hfill
	\includegraphics[width=2.5cm]{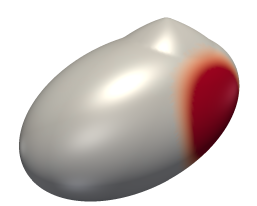}	\hfill
	\includegraphics[width=2.5cm]{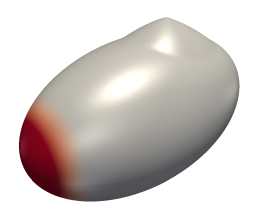}	\hfill
	\includegraphics[width=2.5cm]{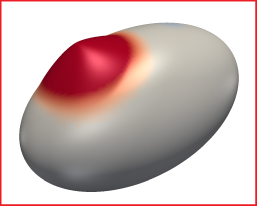}	\hfill
	\includegraphics[width=2.5cm]{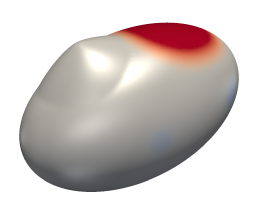}	\\
	\includegraphics[width=2.5cm]{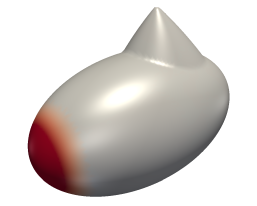}	\hfill
	\includegraphics[width=2.5cm]{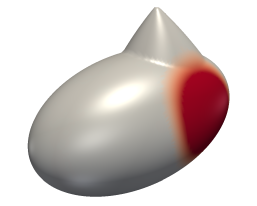}	\hfill
	\includegraphics[width=2.5cm]{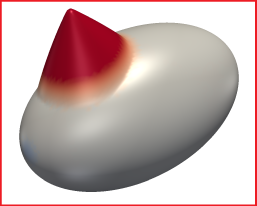}	\hfill
	\includegraphics[width=2.5cm]{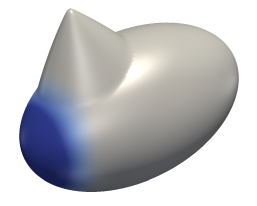}	\hfill
	\includegraphics[width=2.5cm]{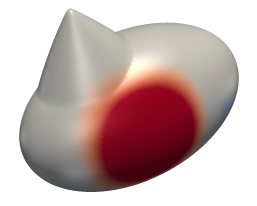}	\\
	
	\caption{Comparison of the first five MHB functions and L$_p$CMs in case of an ellipsoid with a bump. Top row: first 5 non-constant eigenvectors of LBO. Bottom rows: first 5 L$_p$CMs for different bump dimensions.}
	\label{fig:MHB_comp}
\end{figure*}

\tr{In order to further demonstrate how the L$_p$CMs localize the details much better than the Laplacian eigenvectors, 
we consider a synthetic example of an ellipsoid with a growing bump.
The ellipsoid's principal semi-axes are $\{2.5,1.5,1.5\}$ long and it was approximated by
	triangulated mesh of $|V|=16386$ vertices and $|T|=32768$ triangles.
In the top row of Fig. \ref{fig:MHB_comp} we report the first five non-constant eigenvectors
	of LBO corresponding to the first five non-zero eigenvalues obtained by solving the generalized 
	eigenvalue problem \eqref{MHB}. The eigenvectors present global support and neither the first 
	five nor the rest of the eigenvectors, which are not illustrated here for space constraints, are able to localize the bump.
In the bottom rows of Fig. \ref{fig:MHB_comp} we show the first five L$_p$CMs for $p=0.8$ and $\mu=125$,
	for different bump dimensions. In the first and last and row of Fig. \ref{fig:MHB_comp} the bump dimensions correspond.
	%In particular, starting from rows three to five of Figure \ref{fig:MHB_comp} correspond
	%to five L$_p$CMs for an ellipsoid with bump of radius $r=0.2$ and heights $h=\{0.3,0.6,0.9\}$
	%respectively.
%Rows six to eight illustrate the functions computed for a bump of $r=0.6$ and heights $h=\{0.3,0.6,0.9\}$
%and at last, the remaining three rows nine to eleven for radius $r=0.9$ and heights $h=\{0.3,0.6,0.9\}$.
The compact support $L_p$CMs which localize the bumps are highlighted in red boxes. 
}

\smallskip

We conclude this example presenting an empirical investigation on the numerical 
convergence of the proposed ADMM-based minimization scheme.
%
%The ADMM algorithm, in general experiences, shows very slow convergence to high accuracy, moreover, 
%only when the splitting strategy produces sum of convex functions.
%However, in practical applications also a modest level of accuracy is sufficient.

In our formulation (\ref{DP}), we deal with a non-convex orthogonality \tr{constraint} 
and non-convex penalty term, i.e. the sparsity-inducing L$_p$ norm.
Therefore, the convergence to an optimal solution in the global sense is not guaranteed
and we assume that the algorithm converges to at least a local minima,
 which is still a sufficient result for our application.

At that aim, we investigated the empirical convergence via the relative change of the primal variables, 
	and, following \cite{BOYD_ADMM}, the squared primal residual norm $\|\mathbf{r}\|_2^2$ which, 
	according to our implementation, is defined as
	$$ \|\mathbf{r}^{(k)}\|_2^2 = \|\Psi^{(k)}-S^{(k)}\|_F^2 + \|\Psi^{(k)}-E^{(k)}\|_F^2 \,.$$

\begin{figure*}[ht]
	\centering
	\includegraphics[width=7cm]{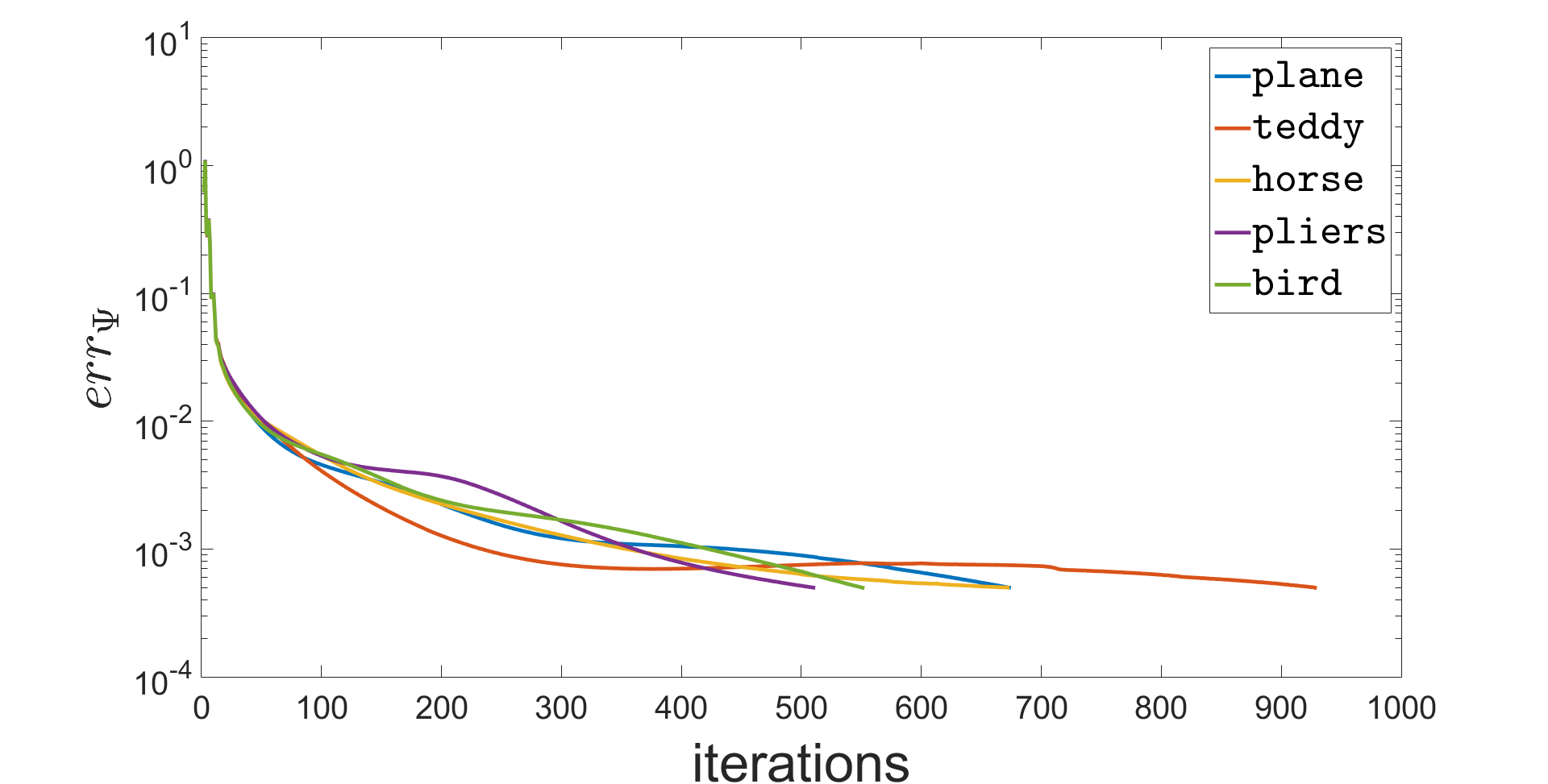}	%\hspace{0.5cm}
	\includegraphics[width=7cm]{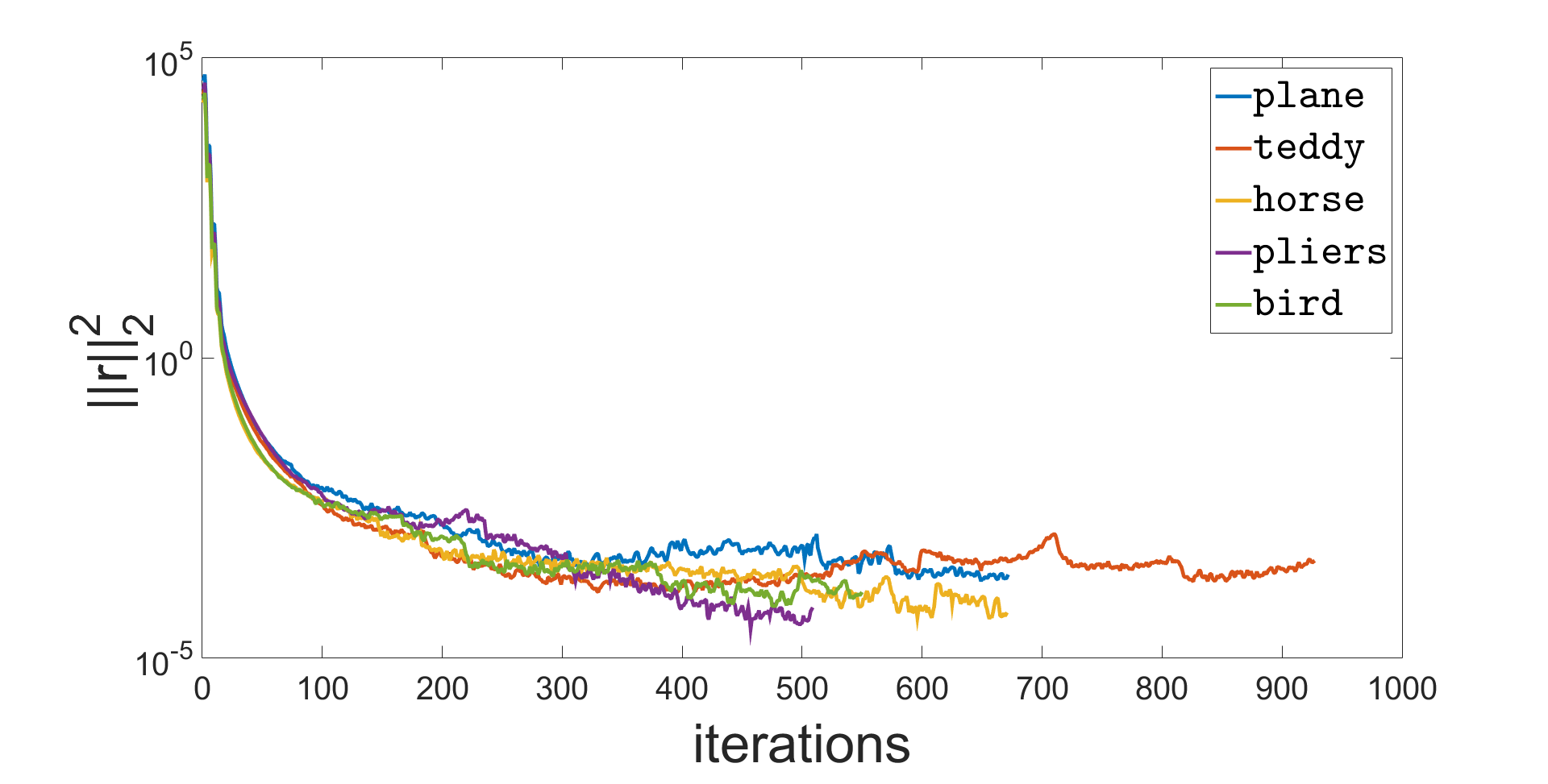}		%\hspace{-0.5cm}
	\caption{ Relative change of primal variable $\Psi$, $err_{\Psi}$, (top) and
		primal residual norm $ \|\mathbf{r}\|_2^2$ (bottom) in terms of ADMM iterations.}
	\label{fig:admm}
\end{figure*}

By the way of illustration, in Figure \ref{fig:admm} we report the convergence plots concerning
some models used for these examples. 
The plots in Fig. \ref{fig:admm}(top)  show that the relative errors $err_{\Psi}$ defined in \eqref{errpsi} on the ADMM iterates $\Psi^{(k)}$ computed
by Step 1 in Algorithm 1,  converge to some limit, which indeed indicates convergence of
the proposed method (at least to local minimizers), whereas the plots in Fig. \ref{fig:admm}(bottom)
demonstrate that the primal residual norms $ \|\mathbf{r}^{(k)}\|_2^2$ reduce. 
%as expected from theory, 
%the approached limits correspond
%to some global/local minimum of the nonconvex nonsmooth TVp functional.

\begin{table}[h]
	\centering
	\caption{\tr{Performance of the mesh segmentation algorithm.}}
	\label{tab:Data}
	\begin{tabular}{lrrrrrr} 
		Data set & $|V|$ & $|T|$ & $K$ & $\mu$	& STEP 1 (s) & STEP 2 (s) 	\\
		\hline 
		\vspace*{-0.3cm} & & & & & & \\
		\texttt{ant}			& 7038	& 14072 & 9	& 150	& 9.69	& 4.79	\\
		\texttt{armadillo}		& 25319	& 50542 & 12& 140	& 47.26	& 13.23	\\
		\texttt{bird}			& 6475	& 12946 & 4	& 300	& 7.32	& 4.54	\\
		\texttt{dolphin}		& 7573	& 15142	& 7 & 150	& 9.40	& 3.48	\\
		\texttt{fawn}			& 3911	& 7818	& 6 & 150	& 4.10	& 2.54	\\
		\texttt{fertility}		& 19994	& 40000	& 7 & 300	& 26.56 & 11.57	\\
		\texttt{fish}			& 5121	& 10238	& 8 & 130	& 7.45	& 2.11	\\
		\texttt{giraffe}		& 9239	& 18474	& 13& 130	& 14.55	& 3.52	\\
		\texttt{glasses}		& 7407  & 14810	& 6 & 150	& 7.75	& 3.36	\\						
		\texttt{hand}			& 6607	& 13210	& 8 & 150	& 8.27	& 3.32	\\
		\texttt{horse}			& 8078	& 16152	& 6 & 300	& 9.59	& 3.81	\\
		\texttt{octopus}		& 5944	& 11888	& 9 & 150	& 11.14	& 4.08	\\
		\texttt{plane}			& 7470	& 14936	& 7 & 150	& 7.54	& 3.34	\\
		\texttt{pliers}			& 3906  & 7808	& 6 & 130	& 5.21	& 2.82	\\
		\texttt{teddy}			& 9548  & 19092	& 7 & 130	& 12.25	& 4.53	\\
		\texttt{teddy\_2}		& 12831 & 25658	& 16& 30	& 22.83	& 4.47	\\
		\texttt{wolf}			& 4712  & 9420	& 7 & 150	& 5.94	& 2.45	\\
	\end{tabular}
\end{table}

\bigskip

\subsection{STEP 2: Mesh Segmentation}

\begin{figure*}[ht]
	\centering
	\includegraphics[width=3.1cm]{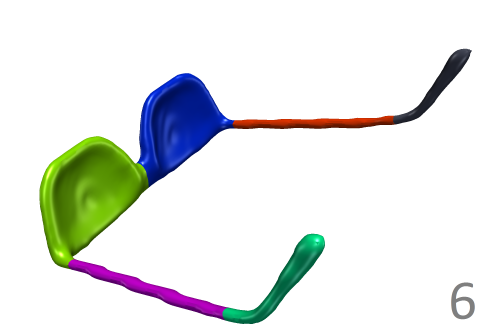}	%\hspace{-0.5cm}
	\includegraphics[width=3.1cm]{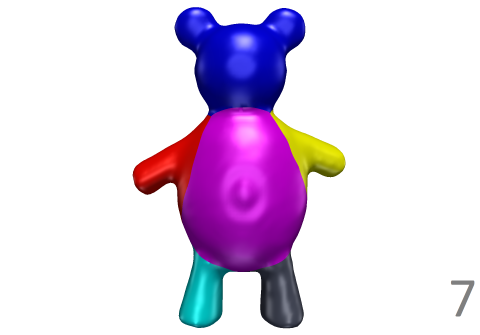}		%\hspace{-0.5cm}
	\includegraphics[width=3.1cm]{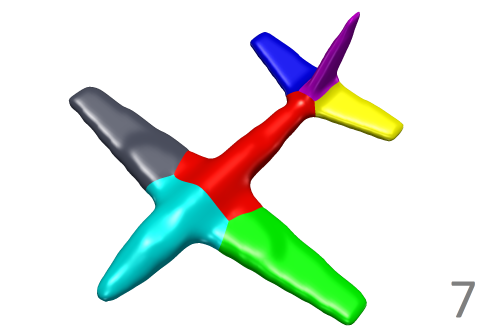}		%\hspace{-0.5cm}
	\includegraphics[width=3.1cm]{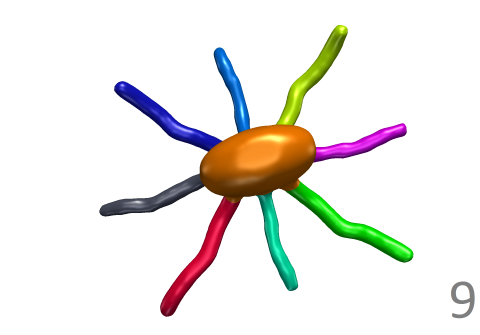}		\\
	\includegraphics[width=3.1cm]{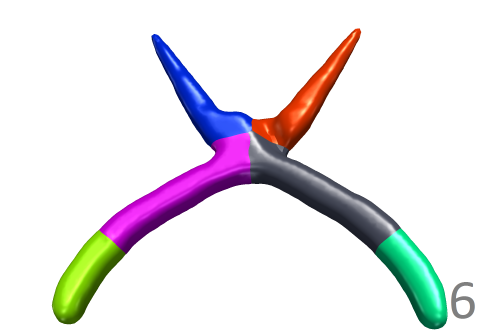}	%\hspace{-0.5cm}
	\includegraphics[width=3.1cm]{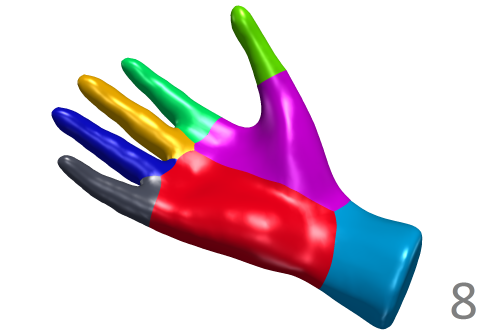}		%\hspace{-0.5cm}
	\includegraphics[width=3.1cm]{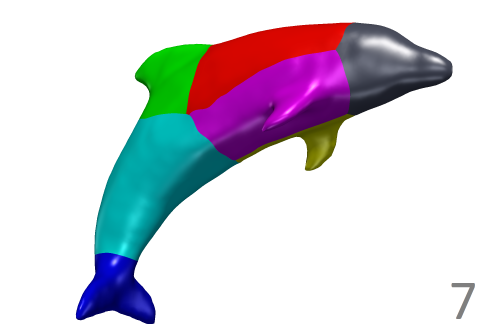}		%\hspace{-0.5cm}
	\includegraphics[width=3.1cm]{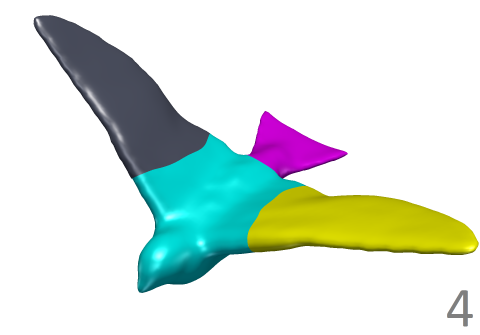}		\\
	\includegraphics[height=2.9cm]{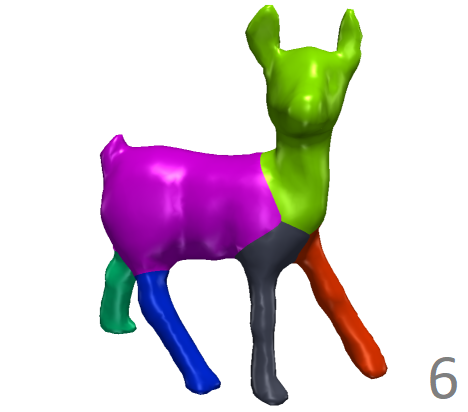}
	\includegraphics[height=2.9cm]{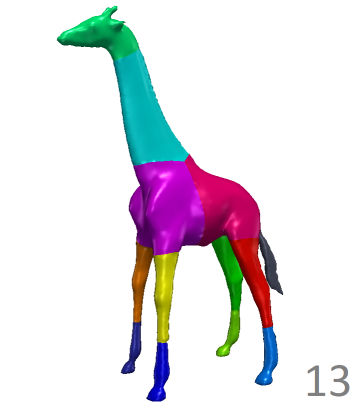}		
	\includegraphics[height=2.9cm]{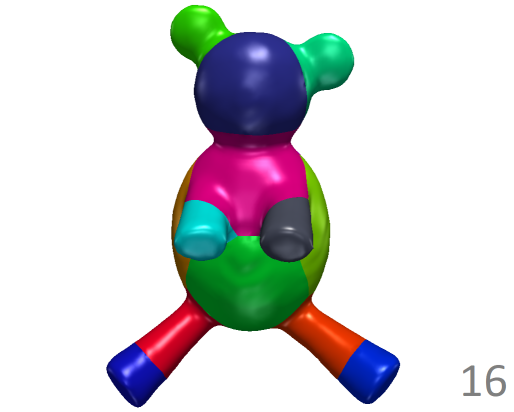}
	\includegraphics[width=3.1cm]{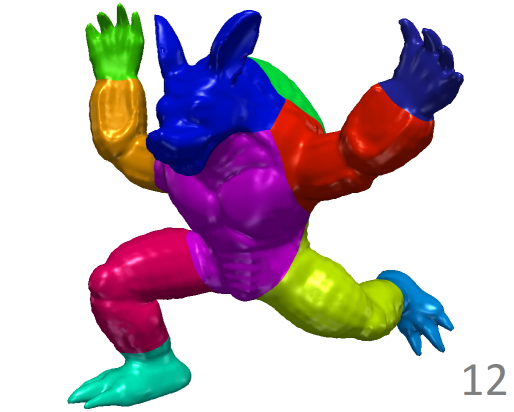}
	
	\caption{Mesh partitioning into salient parts obtained in Step 2 of the \textbf{Algorithm 1}.}
	\label{fig:parts}
\end{figure*}

In Step 2 we apply the region growing algorithm detailed in Section \ref{sec:alg}
to obtain the partition $S=\{S_k\}_{k=1}^N$ which represents a decomposition 
of the mesh into its salient parts.
Several examples of mesh segmentation are shown in Figure \ref{fig:parts}.
The number of partitions produced ($N$) is reported on the bottom right of each segmented object.

\tr{Details of the datasets are given in Table \ref{tab:Data}. In particular, for each mesh, we report the %basis dimension that also corresponds to the 
number of partitions ($N$),
	the value of $\mu$ automatically computed by Step 1b, the time in seconds to obtain the L$_p$CM basis
	of dimension $N$ in Step 1b, and the time for the mesh segmentation procedure in Step 2 of \textbf{Algorithm} 1.} 

It is worth mentioning that our segmentation procedure is naive compared with many other spectral 
segmentation approaches proposed in literature, which are \tr{enriched by many heuristic strategies based on
 curvature  criteria or edge detection, which, however, } can be easily applied also to our basic algorithm.
Nevertheless, the obtained results enhance the good properties of our proposal.
  
\begin{figure*}[ht]
	
	%\hspace*{0.2cm}
	\begin{minipage}{0.7\textwidth}
		\includegraphics[width=1.cm]{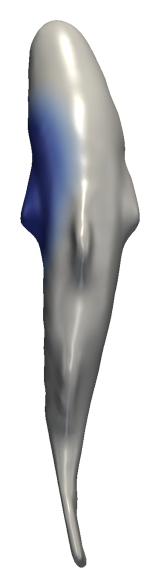}		%\hspace{-0.5cm}
		\includegraphics[width=1.cm]{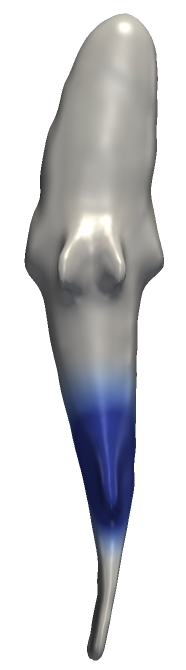}		%\hspace{-0.5cm}
		\includegraphics[width=1.cm]{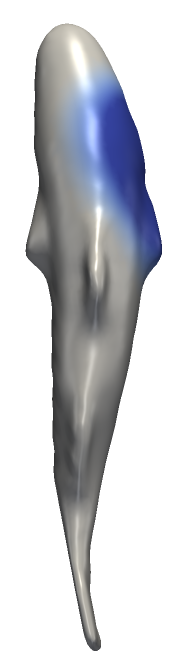}		%\hspace{-0.5cm}
		\includegraphics[width=1.cm]{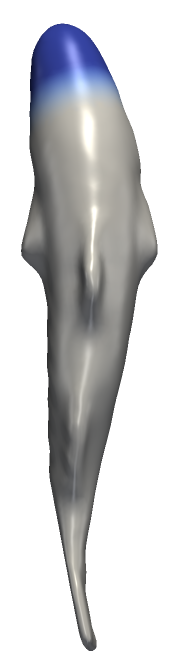}		%\hspace{-0.5cm}
		\includegraphics[width=1.cm]{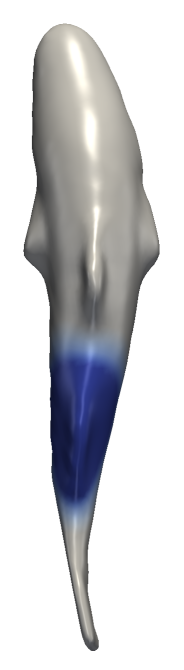}		%\hspace{-0.5cm}
		\includegraphics[width=1.cm]{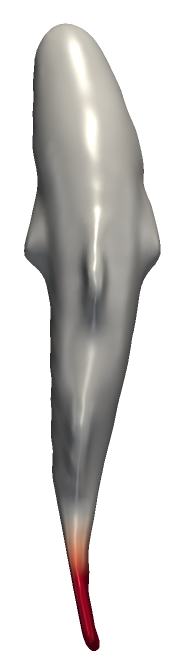}		%\hspace{-0.5cm}
		\includegraphics[width=1.cm]{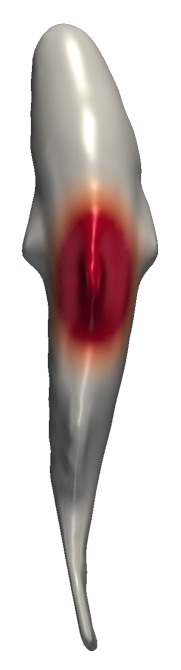}		%\hspace{-0.5cm}
		\includegraphics[width=1.cm]{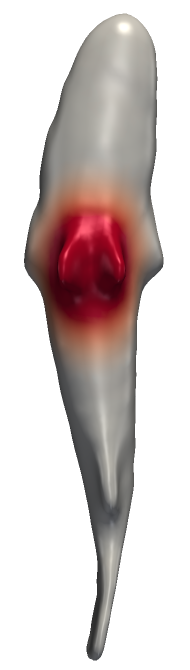}		%\hspace{-0.5cm}
		\hfill \\
		\includegraphics[width=9.1cm]{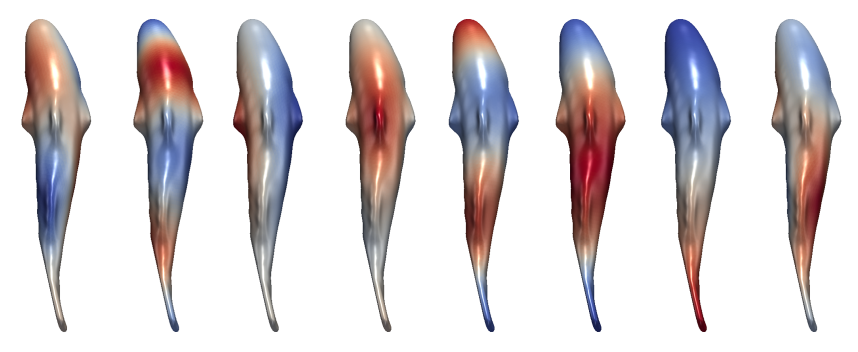}
	\end{minipage}
	\hspace*{0.5cm}
	\vline
	\begin{minipage}{0.2\textwidth}
		\hfill
		\includegraphics[width=1.8cm]{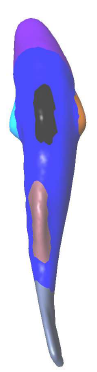}\\		%\hspace{-0.5cm}
	\end{minipage}
	
	\caption{ Supports of the eigenfunctions  for the \texttt{fish} mesh: 
		L$_p$CMs results (top row) and human segmentation (top row, right), eigenfunctions 
		of the affinity matrix proposed in \cite{ZZ} (bottom row).}
	\label{fig:ZZcomp}
\end{figure*}

The model \texttt{fish}, illustrated in Fig.\ref{fig:ZZcomp}, is considered a particularly 
difficult challenge since its featured parts
(fins, head, tail) are smoothly joined with the rest of the body thus
presenting weak boundary strength but good degree of protrusion.
In Fig.\ref{fig:ZZcomp} we show a comparison between our L$_p$CM basis  (top left) 
and the eigenfunctions computed by the truncated spectral decomposition used in \cite{ZZ} (bottom left).
The latter is considered the state-of-the-art among the variational methods using spectral analysis.
 
The salient parts are nicely identified by the L$_p$CMs using only $N=8$ functions, 
mimicking the human driven segmentation shown in Fig.\ref{fig:ZZcomp} (right).
In \cite{ZZ} the authors claim that even for higher space dimensions their method was not able to 
localize the salient parts.
We notice that \tr{in Fig.\ref{fig:ZZcomp} (top row, left) the fish meshes for  $\psi_2$ and $\psi_8$ are visualized upside-down to better} show which fins are localized by these supporting functions.

\begin{figure*}[ht]
	\centering
	\includegraphics[width=5.1cm]{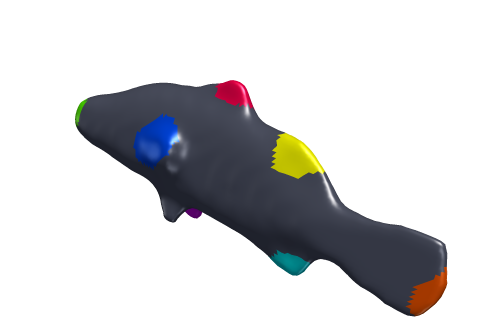}	\hspace{-0.5cm}
	\includegraphics[width=5.1cm]{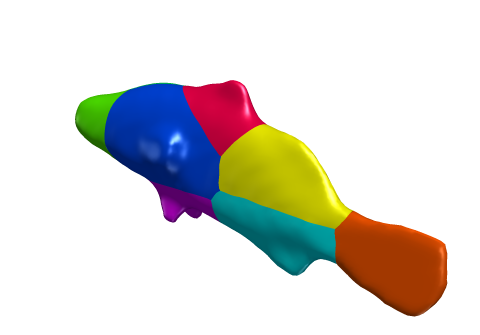}	\hfill
	\includegraphics[width=1.cm]{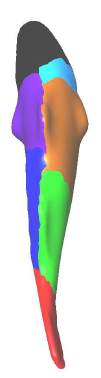}	
	\caption{Segmentation of the \texttt{fish} mesh: seed areas of Step 2(left), mesh segmented by Step 2 (center), and mesh segmentation using \cite{ZZ}.}
	\label{fig:ZZcompRG}
\end{figure*}

 The spectral segmentation results are shown in Fig.\ref{fig:ZZcompRG}.
The starting seeds (left) computed by Step 1 of \textbf{Algorithm 1} are placed correctly,
then the region growing algorithm in Step 2 ends up with the partitioning in Fig.\ref{fig:ZZcompRG}(middle).
%Such result is obtained also due to the fact that the supporting volumes tend to be of similar size, therefore, we are not able to %obtain, in the same basis, all the functions covering just the small fins and a function that covers the rest of the body. 
On the right of Fig.\ref{fig:ZZcompRG} we report the mesh decomposition shown in \cite{ZZ} which has been 
produced with the help of an edge detection strategy introduced in the variational formulation. 
A visual insight allows us to observe some defects for both the top fins,  on the cluster boundaries 
which indeed go through the middle of the fins.

\bigskip

\subsection{STEP 3: Patch-based Partitioning}

The third step of \textbf{Algorithm 1} refines the partitioning obtained by Step 2
finalizing a patch-based manifold partitioning, see Def. 2.
To this end, we first select from $S$ those parts $S_k$ which have genus higher than zero and/or more boundaries, 
and we re-run Step 1 and Step 2 for each of them until every $S_k$ has genus-0 and at most two boundaries.

\begin{figure*}[ht]
	\centering
	\includegraphics[width=3.1cm]{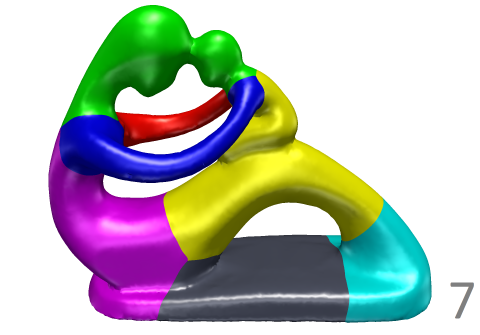}
	\includegraphics[width=3.1cm]{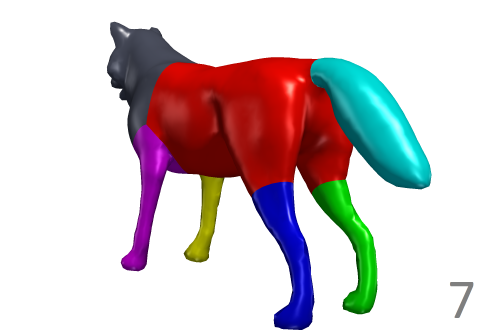}	%\hspace{-0.5cm}
	\includegraphics[width=3.1cm]{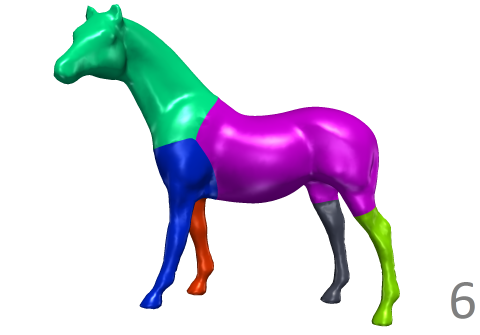}		%\hspace{-0.5cm}
	\includegraphics[width=3.1cm]{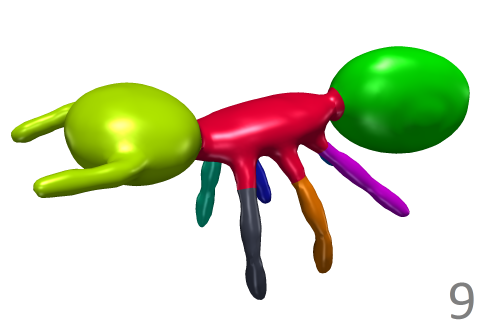}		\\ %\hspace{-0.5cm}
	\includegraphics[width=3.1cm]{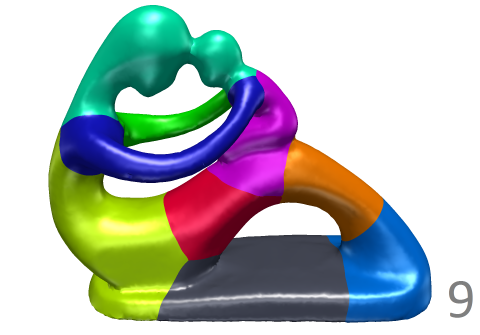}
	\includegraphics[width=3.1cm]{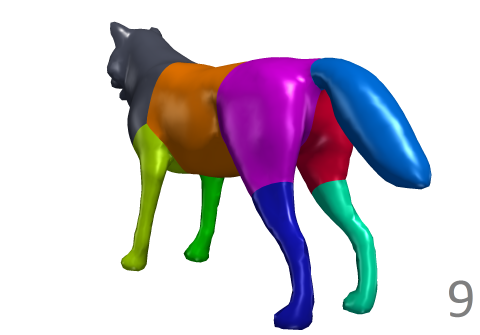}		
	\includegraphics[width=3.1cm]{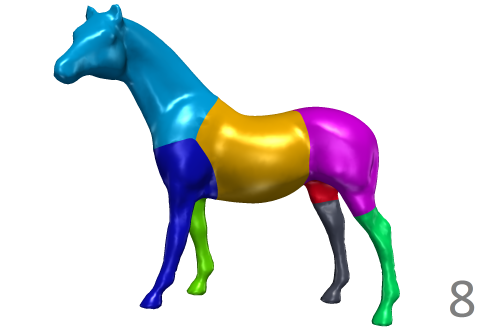}
	\includegraphics[width=3.1cm]{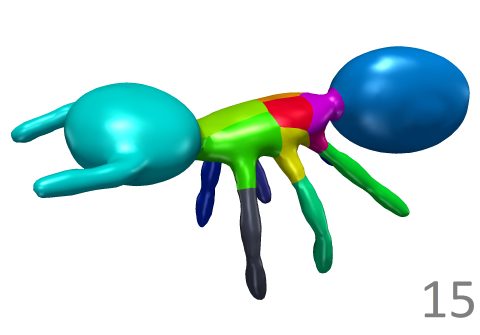}

	\caption{Mesh segmentation into salient parts obtained by applying Step 2 of \textbf{Algorithm 1} (top row);
	patch-based partitioning into genus-0 patches by applying the refinement in Step 3 (bottom row). }
	\label{fig:patches}
\end{figure*}

In Figure \ref{fig:patches} we illustrate a few examples of patch-based partitioning 
resulting from Step 3 (bottom row) in comparison with the mesh partitioning obtained in Step 2 (top row).
For all the meshes reported in this figure, just one part (from left to right yellow/red/magenta/red)
	was further subdivided.
The \texttt{fertility} mesh (left), characterized by four holes, represents a closed mesh of higher genus.
Also in this case, the algorithm was able to both localize salient parts of the mesh and create
	a satisfying genus-0 patching.

\section{Conclusions}
\label{sec:conc}
In this paper, we proposed a sparsity-promoting variational method 
to produce compressed functions L$_p$CMs which are quasi-eigenfunctions of the Laplacian operator.
We proved that the generated functions are highly localized in space and the size of their support
depends on the sparsity parameter $p$ and on the penalty parameter $\mu$. 
%For the proposed variational problem
%which is the key part of the method, 
An Augmented Lagrangian method was applied to solve this non-convex non-differentiable
optimization problem, yielding an iterative algorithm
with efficient solutions to subproblems. 
This compact support basis proves to be very useful for spectral shape processing.  In particular, 
we proposed a unified method for shape partitioning that can be applied to both mesh segmentation and patch-based partitioning, which is based on a more restrictive requirement with respect to mesh segmentation. 
%we applied our method to decompose a manifold represented by a mesh into 
%sub-meshes which can be used either for mesh segmentation or for patch-based partitioning.
In mesh segmentation each sub-mesh represents a meaningful part of the 
object from a human perception point of view. In patch-based partitioning 
instead each part is characterized by genus-0 topology which easily allows it to be
parametrized. 
Our proposed \textbf{Algorithm 1}, thanks to the compactness of the generated L$_p$CM basis, well performs on typical mesh partitioning problems, but it still has some limitations. We have not obtained a convergence proof, which is probably very challenging due to the non-convexity of the problem.
% Like many algorithms, we cannot prove the convergence of the ADMM strategy at present.
Besides, in our tests, we found out that the parameter $p$ 
affects the algorithm efficiency. When $p$ is very small, for example $p < 0.4$, our algorithm
is slow. Acceleration techniques for strong sparsity requirements remain to be designed and they will be 
considered in future work.

\section*{Acknowledgements}
We would like to thank the referees for comments that lead to improvements of the presentation. 
Research was supported in part by the National Group for Scientific Computation (GNCS-INDAM), Research Projects 2015.

\bigskip\bigskip\bigskip\bigskip\bigskip\bigskip
\section*{APPENDIX}

\bigskip

%\\
\textbf{Proof of Theorem\ref{th:cs}.}

\begin{proof}
The proof is decomposed in four steps.
\\
{\bf Step 1: }\\
We first derived the Euler-Lagrange equations for \eqref{problem2}. 
For any fuction $u$, let s(u) denote an element of subdifferential of $\left | u \right |^p$, that is:
\begin{equation}
\label{diffp}
s(u)=sign(u)\cdot  p \cdot \left |u \right |^{p-1}. \; 
\end{equation}
The solutions of (\ref{problem2}) are weak solutions of the following system of nonlinear boundary
value problem:
\begin{equation}
\label{sol}
\frac{1}{\mu} s(\psi_i)-2\lambda_{ii} \psi_i-\Delta \psi_i - \sum_{j \neq i} \lambda_{ij}\psi_j=0, \quad i=1,..,N \; \mbox{on} \quad \Omega
\end{equation}
where $\lambda_{ij}$, with  $\lambda_{ij}=\lambda_{ji}$ are Lagrange multipliers corresponding to orthonormality constraints:
\begin{equation}
\label{OC}
\int_\Omega {\psi_i^2 dx}=1 \quad \mbox{and} \quad \int_\Omega {\psi_i \psi_j dx}=0, \quad \mbox{for}\quad i,j=1,...,N, \; j\ne i.
\end{equation}
\\
{\bf Step 2:} Upper bounds for $\lambda_{ii}$, $\|\psi_i \|_p^p$ and $\|\nabla \psi_i \|_2$\\
For each $i$ multiply both sides of equation (\ref{sol}) by $\psi_i(x)$ and integrate over domain $\Omega$:
\begin{equation}
\int_\Omega{\frac{1}{\mu} s(\psi_i) \psi_i dx}-2\lambda_{ii} \int_\Omega{\psi_i\psi_i dx}-
\int_\Omega{\Delta \psi_i  \psi_i dx} - \sum_{j \neq i} \lambda_{ij}\int_\Omega{\psi_i  \psi_j dx}=0 .
\end{equation}
By using orthonormality conditions (\ref{OC}), we can rewrite the above equation as:
\begin{equation}
\int_\Omega{\frac{1}{\mu} s(\psi_i) \psi_i dx}-2\lambda_{ii} -
\int_\Omega{\Delta \psi_i \psi_i dx}=0 
\end{equation}
that, using integration by parts and zero boundary conditions on $\Omega$, implies that
\begin{equation}
\int_\Omega{\frac{1}{\mu} s(\psi_i)\psi_i dx}-2\lambda_{ii} +
\int_\Omega{\left | \nabla  \psi_i \right| ^2 dx}=0 
\end{equation}
and then
\begin{equation}
\label{lambda}
\lambda_{ii}= \frac{1}{2\mu}  \int_\Omega{s(\psi_i)\psi_i dx}+
\frac{1}{2} \int_\Omega{\left | \nabla  \psi_i \right| ^2 dx}.
\end{equation}
By using definition (\ref{diffp}), relation (\ref{lambda}) can be reformulated as:
\begin{equation}
\label{lambdas}
\lambda_{ii}= \frac{1}{2\mu}  \int_\Omega{p \left |\psi_i \right |^p  dx}+
\frac{1}{2}\int_\Omega{\left | \nabla  \psi_i \right| ^2 dx}
\end{equation}

From Proposition \ref{extProp3_4}, we know that the first compressed mode $\psi$ has support whose
volume satisfy (\ref{volume}). It follows that for $\mu$ sufficiently small and $0<p < 1$, 
the $N$ disjoint copies (i.e. translates) of $\psi$ can be placed in $\Omega$, 
and these $N$ functions are a solution for problem 
(\ref{problem2}). Therefore, in view of Proposition \ref{extProp3_4}, there exist $\mu_0$ (depending on values of $p$,$N$, and $d$) such
that for $\mu < \mu_0$:
\begin{equation}
\sum_{i=1}^N \int_\Omega{ \frac{1}{\mu}\left|{\psi_i}\right|^p dx}+\sum_{i=1}^N \frac{1}{2}\int_\Omega{\left |\nabla {\psi_i} \right |^2 dx} \leq m(\Omega)^{\frac{1}{p}-1} C_1  N \mu^{-\frac{4}{4+d}}.
\end{equation}
Because each of the summands in the left hand side of above inequality is positive, there exist
constant $C_2$ (depending on $d$ and $N$) such that for $\mu < \mu_0$,
\begin{equation}
\label{uppB}
\int_\Omega{ \frac{1}{\mu}\left|{\psi_i}\right|^p dx}\leq m(\Omega)^{\frac{1}{p}-1} C_2  \mu^{-\frac{4}{4+d}} \quad \mbox{and} \quad 
\int_\Omega{\left |\nabla {\psi_i} \right |^2 dx} \leq m(\Omega)^{\frac{1}{p}-1}  C_2  \mu^{-\frac{4}{4+d}}.
\end{equation}

Moreover, replacing the above inequalities into (\ref{lambdas}), it follows that there exist a constant $C_3$
(depending on $d$, $N$ and $p$), such that for $\mu<\mu_0$
\begin{equation}
\label{lambdaiiF}
\left| \lambda_{ii}\right| < p \frac{C_2}{2}  \mu^{-\frac{4}{4+d}}m(\Omega)^{\frac{1}{p}-1}+ p \frac{C_2}{2} \mu^{-\frac{4}{4+d}} m(\Omega)^{\frac{1}{p}-1} <C_3  \mu^{-\frac{4}{4+d}} m(\Omega)^{\frac{1}{p}-1}.
\end{equation}
{\bf Step 3:} Upper bounds for $\lambda_{ij}'s$.\\
Fix $i$. For $k \neq i$, multiply both sides of equation (\ref{sol}) by $\psi_k(x)$ and integrate over $\Omega$:
\begin{equation}
\int_\Omega{\left(\frac{1}{\mu} s(\psi_i)\psi_k -2\lambda_{ii} \psi_i \psi_k-\Delta \psi_i \psi_k -\sum_{j \neq i} \lambda_{ij}\psi_j \psi_k\right) dx}=0
\end{equation}
which, using orthonormality condition (\ref{OC}) and integration by parts, implies that:
\begin{equation}
\frac{1}{\mu} \int_\Omega{s(\psi_i)\psi_k dx}+ \int_\Omega {(\nabla \psi_i) (\nabla \psi_k) dx} -\lambda_{ik}=0.
\end{equation}
Therefore
\begin{equation}
\label{lambdaik}
\lambda_{ik}=\frac{1}{\mu} \int_\Omega{s(\psi_i)\psi_k dx}+ \int_\Omega {(\nabla \psi_i) (\nabla \psi_k) dx}.
\end{equation}
By using relation (\ref{diffp}), we have
\begin{equation}
\left | \frac{1}{\mu} \int_\Omega{s(\psi_i) \psi_k dx} \right|\leq \frac{1}{\mu} \int_\Omega{
\left |s(\psi_i)\psi_k \right|dx}= \frac{p}{\mu} \int_\Omega{
 \left| \psi_i \right|^{p-1}\left|\psi_k\right| dx}= \frac{p}{\mu} \int_\Omega{
 \left| \psi_i \right|^{p} \frac{\left|\psi_k\right|}{\left|\psi_i\right|} dx}.
\end{equation}
%Adesso uso il Primo teorema del valor medio per integrali definiti:	\\
%If $f : [a, b] \rightarrow R$ is continuous and $g$ is an integrable function that does not change sign on $[a, b]$, then there exists $c \in (a, b)$ such that
%$\int_a^b{   f ( x ) g ( x )  d x} = f ( c ) \int_a^b {  g ( x )  d x }.$
\\
Since $\left|\psi_i\right|^{p}$ does not change sign on $\Omega$, by the First Mean Value Theorem for Integrals,
there exists $\xi \in \Omega$, with $\psi_i(\xi) \neq 0$, such that, if we set
$M=\frac{\left|\psi_k(\xi)\right|}{\left|\psi_i(\xi)\right|}$, it follows that
$$  \frac{p}{\mu} \int_\Omega{
 \left| \psi_i \right|^{p} \frac{\left|\psi_k\right|}{\left|\psi_i\right|} dx}=M  \frac{p}{\mu} \int_\Omega{
 \left|\psi_i \right|^{p} dx}.
$$

 Making use of (\ref{uppB}), we conclude that

\begin{equation}
\label{eq:up1}
\left| \frac{1}{\mu} \int_\Omega{s(\psi_i)\psi_k dx} \right| \leq p  M  m(\Omega)^{\frac{1}{p}-1} C_2   \mu^{-\frac{4}{4+d}}.
\end{equation}

Finally, using Cauchy-Schwarz and equation (\ref{uppB}),
\begin{multline}
\label{eq:up2}
\left | \int_\Omega{(\nabla \psi_i)  (\nabla \psi_k) dx} \right | \leq \left(
\int_\Omega{\left| \nabla \psi_i \right|^2 dx} \right)^{\frac{1}{2}} \left(
\int_\Omega{\left| \nabla \psi_k \right|^2 dx} \right)^{\frac{1}{2}} < \\(m(\Omega)^{\frac{1}{p}-1}  C_2  \mu^{-\frac{4}{4+d}})^{\frac{1}{2}} (m(\Omega)^{\frac{1}{p}-1}  C_2  \mu^{-\frac{4}{4+d}})^{\frac{1}{2}} \, = \, m(\Omega)^{\frac{1}{p}-1}  C_2  \mu^{-\frac{4}{4+d}}.
\end{multline}
\\
Substituting the two upper bounds given in \eqref{eq:up1} and \eqref{eq:up2} into equation (\ref{lambdaik}), we have 
for $\mu < \mu_0$
\begin{multline}
\label{lambdaikF}
\left | \lambda_{ik} \right | < p M   m(\Omega)^{\frac{1}{p}-1} C_2   \mu^{-\frac{4}{4+d}} + m(\Omega)^{\frac{1}{\tilde{p}}-1} C_2 \mu^{-\frac{4}{4+d}}=\\
 m(\Omega)^{\frac{1}{p}-1} C_2  \mu^{-\frac{4}{4+d}}(p M+1)< C_4 m(\Omega)^{\frac{1}{\tilde{p}}-1}   \mu^{-\frac{4}{4+d}}
\end{multline}
where $C_4$ depends on $N$, $M$, $p$ and $\mu$.

\smallskip

\noindent{\bf Step 4:} Bounding the volume of the support $\psi_i$'s\\
For each $i$ multiply both sides of equation (\ref{sol}) by
$\frac{1}{p} sign(\psi_i) \left| \psi_i \right|^{1-p}$ and integrate over domain $\Omega$:
\begin{equation}
\label{support}
\frac{1}{\mu} \left| supp(\psi_i) \right|-\frac{2}{p}\lambda_{ii} \int_\Omega{\left | \psi_i \right|^{2-p}dx}-   \frac{1}{p}\int_\Omega{\Delta \psi_i sign(\psi_i) \left|\psi_i\right|^{1-p}dx} -\frac{1}{p}\sum_{j\neq i} \lambda_{ij} \int_\Omega{\psi_j sign(\psi_i) \left|\psi_i\right|^{1-p}dx}=0,
\end{equation}
namely,
\begin{multline}
\label{support1}
\frac{1}{\mu} \left| supp(\psi_i) \right|=\\ \left|\frac{2}{p}\lambda_{ii} \int_\Omega{\left | \psi_i \right|^{2-p}dx}+   \frac{1}{p}\int_\Omega{\Delta \psi_i  sign(\psi_i) \left|\psi_i\right|^{1-p}dx}+\frac{1}{p}\sum_{j\neq i} \lambda_{ij} \int_\Omega{\psi_j  sign(\psi_i)  \left|\psi_i\right|^{1-p}dx}\right|
\leq\\ \left|\frac{2}{p}\lambda_{ii} \int_\Omega{\left | \psi_i \right|^{2-p}dx}+   \frac{1}{p}\int_\Omega{\Delta \psi_i sign(\psi_i) \left|\psi_i\right|^{1-p}dx}\right|+
\left|\frac{1}{p}\sum_{j\neq i} \lambda_{ij} \int_\Omega{\psi_j  sign(\psi_i) \left|\psi_i\right|^{1-p}dx}\right|.
\end{multline}
Define 
$$\Omega^+=\left \{x \in \Omega : \psi_i(x)>0\right\}$$
and 
$$\Omega^-=\left \{x \in \Omega : \psi_i(x)<0\right\}.$$
According to Green's formula
$$\int_{\Omega^+}{\Delta \psi_i dx}=\int_{\partial \Omega^+}{\frac{\partial\psi_i}{\partial \nu} dS} \leq 0,$$
where $\nu$ is outward pointing unit normal vector along $\partial \Omega^+$. Since $\psi$ is positive
in $\Omega^+$ and becomes zero on $\partial \Omega^+$, the right-hand side of the above expression is not positive. \\
With a similar argument, we have that
$$\int_{\Omega^-}{\Delta \psi_i dx}=\int_{\partial \Omega^-}{\frac{\partial\psi_i}{\partial \nu} dS} \geq 0.$$
Hence, since $\left| \psi_i\right|^{1-p} \geq 0$ $\forall i$, it follows that:
\begin{equation}
\label{green}
\int_\Omega{\Delta \psi_i sign(\psi_i) \left|\psi_i\right|^{1-p}dx}=\int_{\Omega^+}{\Delta \psi_i \left|\psi_i\right|^{1-p}dx}-\int_{\Omega^-}{\Delta \psi_i \left|\psi_i\right|^{1-p}dx}\leq 0.
\end{equation}

Using inequality (\ref{green}), (\ref{support1}) can be rewritten as:
\begin{multline}
\label{support21}
\frac{1}{\mu} \left| supp(\psi_i) \right|\leq \left|\frac{2}{p}\lambda_{ii} \int_\Omega{\left | \psi_i \right|^{2-p}dx}\right|+ \left|\frac{1}{p}\sum_{j\neq i} \lambda_{ij} \int_\Omega{\left |\psi_j \right|\left|\psi_i\right|^{1-p}dx}\right|\leq \\
\frac{2}{p}|\lambda_{ii}| \int_\Omega{\left | \psi_i \right|^{2-p}dx}+ \frac{1}{p}\sum_{j\neq i} |\lambda_{ij}| \int_\Omega{\left |\psi_j \right| \left|\psi_i\right|^{1-p}dx} \leq 
\\ \frac{2}{p} |\lambda_{ii}| \int_\Omega{\left | \psi_i \right|\left| \psi_i \right|^{1-p}dx}+ \frac{1}{p}\sum_{j\neq i} |\lambda_{ij}| \int_\Omega{\left |\psi_j \right| \left|\psi_i\right|^{1-p}dx}
\end{multline}
We set $\tilde{p}=1-p$,  $0<\tilde{p}<1$ , for $0<p<1$.
\begin{equation}
\label{support2}
\frac{1}{\mu} \left| supp(\psi_i) \right|\leq \frac{2}{p} |\lambda_{ii}| \int_\Omega{\left | \psi_i \right|\left| \psi_i \right|^{\tilde{p}}dx}+ \frac{1}{p}\sum_{j\neq i} |\lambda_{ij}| \int_\Omega{\left |\psi_j \right| \left|\psi_i\right|^{\tilde{p}}dx}
\end{equation}
Since $\left|\psi_i\right|^{\tilde{p}}$  does not change sign on $\Omega$, by the First Mean Value Theorem for Integrals, there exist $\xi, \eta \in \Omega$ such that, if we set $\bar{M}=\left|\psi_i(\xi)\right|$ and
 $\tilde{M}=\left|\psi_j(\eta)\right|$, relation (\ref{support2}) can be rewritten as:
\begin{equation}
\label{support3}
\frac{1}{\mu} \left| supp(\psi_i) \right|\leq \frac{2}{p}  \bar{M} |\lambda_{ii}|  \int_\Omega{ \left| \psi_i \right|^{\tilde{p}}dx}+ \frac{1}{p}\sum_{j\neq i}  \tilde{M} |\lambda_{ij}| \int_\Omega{\left|\psi_i\right|^{\tilde{p}}dx}.
\end{equation}
By using (\ref{uppB}), (\ref{lambdaiiF}),(\ref{lambdaikF}), then (\ref{support3}) becomes:
\begin{multline}
\label{support4}
\frac{1}{\mu} \left| supp(\psi_i) \right|\leq \frac{2}{p}  \bar{M} C_3 m(\Omega)^{\frac{1}{p}-1} \mu^{-\frac{4}{4+d}} 
m(\Omega)^{\frac{1}{1-p}-1}C_2 \mu^{-\frac{4}{4+d}+1} +\\ \frac{1}{p} (N-1)  \tilde{M} 
C_4 m(\Omega)^{\frac{1}{p}-1}\mu^{-\frac{4}{4+d}} m(\Omega)^{\frac{1}{1-p}-1}C_2 \mu^{-\frac{4}{4+d}+1}
\\ \leq C_5 \mu^{-\frac{8}{4+d}+1} m(\Omega) ^{\frac{1}{p(1-p)}-2}
\end{multline}
where $C_5$ depends on $N$ and $p$.
\end{proof}
\end{document}